\documentclass[journal]{IEEEtran}
\usepackage{amsfonts}
\usepackage{amsmath}
\usepackage{amssymb}
\usepackage{xcolor}
\usepackage{graphicx}
\usepackage{cite}
\usepackage{array}
\usepackage{subfigure}
\usepackage{booktabs}
\usepackage{enumerate}
\usepackage{color}
\usepackage{algpseudocode}
\usepackage{epstopdf}
\usepackage{ntheorem}
\usepackage{pifont}
\usepackage{setspace}
\usepackage{url}
\ifCLASSINFOpdf
\else
\fi
\newtheorem{theorem}{Theorem}

\newtheorem{assumption}{Assumption}

\newtheorem{corollary}{Corollary}

\newtheorem{lemma}{Lemma}

\newtheorem{problem}{Problem}

\newtheorem{remark}{Remark}

\newtheorem{definition}{Definition}
\newenvironment{proof}[1][Proof]{\noindent \textbf{#1.} }{\  \rule{0.5em}{0.5em}}
\ifCLASSINFOpdf \else \fi
\allowdisplaybreaks[4]

\begin{document}
\begin{spacing}{1.0}
\title{Using delays for prescribed-time static output feedback stabilization of LTI controllable and observable systems}
\author{Kang-Kang Zhang and Emilia Fridman, \IEEEmembership{Fellow, IEEE}
\thanks{This work was supported in part by the Azrieli International Postdoctoral Fellowship, in part by the Israel Science Foundation
under Grant 446/24, and in part by the Chana and Heinrich Manderman
Chair at Tel Aviv University. \textit{(Corresponding author: Kang-Kang Zhang)}}
\thanks{The authors are with the School of Electrical and Computer Engineering, Tel Aviv University, 6997801, Israel (e-mail: kangkang\_kkz@163.com; emilia@tauex.tau.ac.il)}}

\maketitle

\begin{abstract}
This paper addresses a stabilization problem:
can every controllable and observable LTI system be
stabilized by static output feedback (SOF)?
To the best of our knowledge, we provide, for the first
time, a solution to this problem.
Furthermore, the proposed design achieves prescribed-time
convergence to the origin independently of
the initial condition.
Two SOF methods are developed.
The first employs multiple artificial discrete output delays
with smooth time-varying gains and provides an explicit
gain construction for every controllable and observable
LTI system, based on a generalized
controllability--observability Gramian introduced
in this paper.
The second combines constant-gain feedback of the current
output with time-varying-gain feedback of a single
artificially delayed output and also applies to every
controllable and observable LTI system.
Its simpler structure requires auxiliary parameters
whose existence is guaranteed but whose explicit
construction is provided only for SISO systems.
The multiple-delay design is further extended to known
bounded and fast-varying (without any constraints on the delay derivative) measurement delays that may be large,
yielding fixed-time stabilization under a corresponding
rank condition, with an assignable settling-time bound
greater than three times the delay bound.
For practical implementation, we further present an SOF controller with time-varying gains that depends on discrete-time delayed measurements. For unknown constant measurement or input delays, fixed-time
stabilization is established for systems with nilpotent
state matrices.
A numerical implementation of Gramian-based gains
is also suggested.
Numerical examples involving nonminimum-phase,
flexible-structure, and cart--pole systems illustrate
the effectiveness of the proposed methods.
\end{abstract}


\section{Introduction}

Static output feedback (SOF) for LTI systems is attractive because of its simple
implementation and low computational cost, but its stabilizing
capability is generally more limited than that of state
feedback \cite{Syrmos97auto,sadabadi16}. These limitations have motivated
extensions of SOF through time-varying
gains and artificial delays. For the constant-gain case, recent results
\cite{lofberg26,barkai26,ahmadi26} established that deciding
whether a linear system admits a stabilizing constant
output-feedback gain is NP-hard.
\subsection{Time-Varying Gains and Artificial Delays for SOF}
Consider the LTI system
\begin{equation}
\left \{
\begin{array}
[c]{rl}%
\dot{x}(t) & =Ax(t)+Bu(t),\\
y(t) & =Cx(t),
\end{array}
\right.  \label{eq:MIMO_plant}%
\end{equation}
where $x(t)\in \mathbb{R}^{n}$, $u(t)\in \mathbb{R}^{m}$, $y(t)\in \mathbb{R}%
^{p}$, $A\in \mathbb{R}^{n\times n}$, $B\in \mathbb{R}^{n\times m}$, and
$C\in \mathbb{R}^{p\times n}$.

Brockett posed the following
time-varying output-feedback stabilization problem
\cite{Brockett99}.

\medskip
\noindent
\textbf{Brockett problem.}
Given matrices $A$, $B$, and $C$, under what circumstances
does there exist a time-varying controller
\begin{equation}
u(t)=K(t)y(t)
\label{eq:time_varying_output_feedback}
\end{equation}
such that the closed-loop system (\ref{eq:MIMO_plant}) and (\ref{eq:time_varying_output_feedback}) is asymptotically stable?
\medskip

This problem motivated the development of periodic
output-feedback methods
\cite{Moreau99cdc,Moreau04scl,Allwright05auto,Michiels04ima}. In this direction, \cite{boikov05arc} presented stabilization
conditions and constructive methods using delay-free
time-varying output feedback for a broad class of systems.
 For SISO systems, memoryless output feedback $u(t)=-k(t)y(t)$ with a time-varying gain $k(t)$
and delayed output feedback $u(t)=-ky(t-\tau)$ with constant $k$ and $\tau>0$ have complementary stabilizing capabilities \cite{Michiels04ima}. The latter paper provided a complete comparison of two feedbacks for second-order SISO systems.
For systems with a given feedback delay $\tau>0$, the Brockett
problem was extended in \cite{Insperger08ifac} as follows.

\medskip
\noindent
\textbf{Brockett problem with feedback delay.}
Given matrices $A$, $B$, and $C$ and a feedback delay
$\tau>0$, under what circumstances does there exist a
time-varying controller
\begin{equation}
u(t)=K(t)y(t-\tau)
\label{eq:time_varying_delayed_output_feedback}
\end{equation}
such that the closed-loop system
(\ref{eq:MIMO_plant}) and
(\ref{eq:time_varying_delayed_output_feedback})
is asymptotically stable?
\medskip

To address this problem and the difficulty caused by
the infinite-dimensional closed-loop dynamics,
\cite{Insperger08ifac} investigated the act-and-wait
control framework introduced in \cite{Insperger06tcst},
in which the feedback is periodically switched off
and on. When the waiting interval is at least as long
as the feedback delay, stability can be characterized
by a finite-dimensional monodromy matrix.
However, general conditions for the existence of a
stabilizing gain within this framework were left
as an open problem.

Multiple artificial delays have also been investigated
for SOF stabilization of special classes
of linear systems using frequency-domain methods
\cite{Gomez19auto,nicule04tac,kharitonov05tac,michiels2007}
and time-domain methods
\cite{karafyllis08scl,emiliabook,fridman25delaycontrol,emilia16auto,emilia17auto,anton18auto,efimov20automatica,iftime26}.
These results primarily establish asymptotic, exponential,
or hyperexponential stability under restrictions
on the system structure.

Recently, artificial time delays have been incorporated into
smooth feedback via the act-and-wait method to achieve
convergence within an initial-condition-independent time
bound, yielding prescribed-time stabilization for systems
without input delays and fixed-time stabilization for
systems with input delays \cite{zhoutac}.
Using this method, \cite{zhang26tac} studied fixed-time
stabilization of systems with nilpotent state matrices
subject to large unknown delays. In addition, an alternative approach to fixed-time
stabilization based on sampled-data control was proposed
in \cite{Mazenc26fixed}.

Prescribed-time SOF design is particularly challenging because
the full state information is unavailable. By employing
Taylor expansions or finite-difference approximations,
in our recent paper \cite{zhang26auto} we designed such a prescribed-time SOF controller for LTI systems. However, its applicability is restricted
to systems without internal dynamics \cite{Zhou23automatica}, where the observability assumption is not sufficient (e.g. not applicable to nonminimum-phase systems).
A systematic construction of SOF controllers for
stabilization of general controllable and observable linear systems
remains largely unexplored.

\subsection{Problem Formulation and Contributions}
In this paper, to the best of our knowledge, we for the first time solve the stabilization problem for \emph{every} controllable and observable LTI system
(\ref{eq:MIMO_plant}) by SOF using artificial delays. Beyond stabilization, the proposed design guarantees
exact convergence to the origin within any prescribed time, independently of the initial condition.

Following \cite{polyakov12,zhou2021tac,song2017auto,ye2023tac,zhoutac}, we adopt the following
definition of fixed-time and prescribed-time stability.
\begin{definition}
\label{deff}
Consider a closed-loop system with state
$x(t;x_0)\in\mathbb R^n$.
The origin is said to be fixed-time stable if it is
Lyapunov stable and there exists a constant $T>0$,
independent of $x_0$, such that
$x(t;x_0)=0$ for all $t\geq T$ and every
$x_0\in\mathbb R^n$.
The constant $T$ is called a settling-time bound.
If the settling-time bound $T$ can be arbitrarily
prescribed a priori through the controller parameters,
then the origin is said to be prescribed-time stable.
\end{definition}

The following assumption is imposed.
\begin{assumption}
\label{assumption_minimal} The pair $(A,B)$ is controllable, and the pair
$(A,C)$ is observable.
\end{assumption}
The prescribed-time stabilization problem is formulated
as follows.
\begin{problem}
Consider system (\ref{eq:MIMO_plant}) under Assumption \ref{assumption_minimal}. Given a prescribed time $T>0$, the objective is to design an
artificial-delay-based time-varying SOF controller
such
that the closed-loop system is prescribed-time stable in the sense of
Definition \ref{deff}.
\end{problem}

In this paper, we propose two SOF
methods for prescribed-time stabilization of linear systems.
The first employs multiple artificial output delays
and provides an explicit gain construction for
\emph{every} controllable and observable LTI system,
based on a generalized controllability--observability
Gramian introduced in this paper.
The second combines a constant-gain feedback of the current
output with a time-varying-gain feedback of a single
artificially delayed output and also applies to
\emph{every} controllable and observable LTI system.
It has a simpler structure but requires auxiliary
parameters whose existence is guaranteed, while their explicit
construction is provided only for SISO systems.
The multiple-delay design is further extended to known
bounded and fast-varying (without any constraints on the delay derivative)
measurement delays that may be large, yielding fixed-time stabilization
under a corresponding rank condition.
For unknown constant measurement or input delays, fixed-time
stabilization is established for systems with nilpotent
state matrices, thereby extending the state-feedback
results of \cite{zhang26tac} to output-feedback control.

The remainder of this paper is organized as follows.
Section~\ref{multiple_output} presents the multiple-delay SOF design.
Section~\ref{Ins_output} presents an SOF controller that uses
the current output and a single delayed output. Section~\ref{sec:measurement_delays} extends the multiple-delay design
to systems with measurement delays.
Section~\ref{subsec:finite_dimensional_implementation}
presents a numerical implementation of Gramian-based gains. Section~\ref{Num_ex} provides numerical examples, and Section \ref{Conclu}
concludes the paper. The proofs of auxiliary results
are collected in the Appendix.

\textbf{Notation.} $\mathbb{Z}_{0}^{+}$ and $\mathbb{N}$ denote the sets of
nonnegative and positive integers, respectively. For a vector or matrix, $\|\cdot\|$ denotes the Euclidean norm or
its induced matrix norm, respectively. For a matrix $M$, $M^{\mathrm T}$,
$\operatorname{rank}(M)$, $\det(M)$,
$\operatorname{adj}(M)$, $\lambda_{\min}(M)$,
$\lambda_{\max}(M)$, and $\sigma_{\min}(M)$ denote
its transpose, rank, determinant, adjugate,
minimum eigenvalue, maximum eigenvalue, and
minimum singular value, respectively. The symbols $\operatorname{Im}(M)$ and $\operatorname{span}%
\{ \cdot \}$ denote the image of $M$ and the linear span of a set, respectively.
The operators
$\operatorname{vec}(\cdot)$ and $\operatorname{vec}^{-1}(\cdot)$ denote the
column-wise vectorization of a matrix and its inverse operation, respectively.
Moreover, $\lfloor a\rfloor$ denotes the greatest
integer not exceeding $a$. For an interval \(\mathcal I\subset\mathbb R\) and a nonnegative integer
\(\nu\), \(\mathbf{C}^{\nu}(\mathcal I)\) denotes the space of
\(\nu\)-times continuously differentiable real-valued functions on
\(\mathcal I\). The operator $\operatorname{col}(\cdot)$ denotes vertical
concatenation of its arguments, and $\otimes$ denotes
the Kronecker product.

%

\section{Feedback with Multiple Output Delays}
\label{multiple_output}
Given a prescribed time $T>0$, let $h=T/3$. Define the cycle index $l$ and the
local phase $\theta$ by
\begin{equation}
l=\left\lfloor \frac{t}{3h}\right\rfloor,
\qquad
\theta=t-3hl\in[0,3h).
\label{eq:cycle_index_phase}
\end{equation}
The act-and-wait function is defined as
\begin{equation}
g(t)=
\begin{cases}
0, & \theta\in[0,2h],\\
1, & \theta\in(2h,3h).
\end{cases}
\label{eq:act_and_wait_function}
\end{equation}
Throughout the paper, each controller is defined to give $u(t)=0$
whenever $g(t)=0$, without evaluating the delayed measurements
or gain functions outside their defined intervals.

Choose an integer $n_0$ with $1\leq n_0\leq n$
and artificial delays $h_1,\ldots,h_{n_0}$ satisfying
$h\leq h_1<h_2<\cdots<h_{n_0}\leq2h$. Define the delayed observability matrix
\begin{equation}
\mathcal{O}_{\boldsymbol{h}}
=
\begin{bmatrix}
C\mathrm{e}^{-Ah_{1}}\\
C\mathrm{e}^{-Ah_{2}}\\
\vdots\\
C\mathrm{e}^{-Ah_{n_{0}}}
\end{bmatrix}
\in\mathbb{R}^{n_{0}p\times n}.
\label{eq:sampled_observability_matrix}
\end{equation}
For SISO systems, one typically chooses $n_{0}=n$.
For MIMO systems with $p>1$, fewer delays may suffice,
allowing $n_0<n$, as illustrated later in
Section~\ref{subsec:MIMO_nilpotent_example}.
\begin{lemma}
\label{cor_general_mimo}
Suppose that $(A,C)$ is observable.
Then there exist an integer $n_0\in\{1,\ldots,n\}$
and artificial delays
$h\leq h_1<\cdots<h_{n_0}\leq2h$ such that
\begin{equation}
\operatorname{rank}\!\left(\mathcal O_{\boldsymbol h}\right)
=n.
\label{eq:sampled_observability_rank}
\end{equation}
\end{lemma}

\begin{proof}
Let $\mathcal I=[h,2h]$. We first show that
\[
\operatorname{span}\left\{
\operatorname{Im}\left(
\mathrm e^{-A^{\mathrm T}s}C^{\mathrm T}
\right):s\in\mathcal I
\right\}
=\mathbb R^n.
\]
Otherwise, there would exist a nonzero vector
$v\in\mathbb R^n$ such that
$C\mathrm e^{-As}v=0$ for every $s\in\mathcal I$.
Differentiating at any $s_0\in(h,2h)$ gives
$CA^k\mathrm e^{-As_0}v=0$ for $k=0,\ldots,n-1$.
Observability of $(A,C)$ implies
$\mathrm e^{-As_0}v=0$.
Since the matrix exponential is nonsingular,
this contradicts $v\neq0$.

From this spanning family, select $n$ linearly
independent vectors. These vectors belong to
$\operatorname{Im}(\mathrm e^{-A^{\mathrm T}s}C^{\mathrm T})$
at at most $n$ distinct values of $s$.
Denote these values, arranged in increasing order,
by $h_1,\ldots,h_{n_0}$, where $1\leq n_0\leq n$.
The columns of $\mathcal O_{\boldsymbol h}^{\mathrm T}$
span all the selected vectors, so
$\operatorname{rank}\mathcal O_{\boldsymbol h}=n$.
\end{proof}

Next, we provide a constructive method for selecting $n$ artificial
delays $h_i$, $i=1,\ldots,n$, such that
(\ref{eq:sampled_observability_rank}) holds.
\begin{lemma}
\label{cor:equally_spaced_delays}
Suppose that $(A,C)$ is observable. Let $h_{0}\in[h,2h)$ and choose
$h_{i}=h_{0}+(i-1)d$, $i=1,\ldots,n$,
with some $d>0$ that satisfies $h_{0}+(n-1)d<2h$.
Then there exists a sufficiently small constant $d^{\ast}>0$ such that, for every
$0<d<d^{\ast}$, (\ref{eq:sampled_observability_rank}) holds with $n_{0}=n$.
\end{lemma}

See Appendix A1 for proof.
\begin{remark}
As shown in Appendix A1, for $n\geq2$, $d^\ast$ can be
chosen according to (\ref{eq:admissible_delay_spacing}).
For $n=1$, only one artificial delay is needed,
so no spacing condition is required.
\end{remark}
For $s\in[2h,3h]$, define
\begin{equation}
\mathit{\Gamma}_{\mathrm d,i}(s)
:=
\mathrm e^{A^{\mathrm T}(s-h_i-2h)}C^{\mathrm T}
\otimes \mathrm e^{A(3h-s)}B
\in\mathbb R^{n^2\times mp},
\label{eq:multiple_delay_lifted_block}
\end{equation}
where $i=1,\ldots,n_{0}$, and collect these block matrices as
\begin{equation}
\mathit{\Gamma}_{\mathrm{d}}(s)
:=
\begin{bmatrix}
\mathit{\Gamma}_{\mathrm{d},1}(s)
&
\cdots
&
\mathit{\Gamma}_{\mathrm{d},n_{0}}(s)
\end{bmatrix}
\in\mathbb{R}^{n^{2}\times n_{0}mp}.
\label{eq:multiple_delay_lifted_matrix}
\end{equation}

For a positive integer $\nu\geq1$, let $S$ satisfy
\begin{equation}
\begin{aligned}
\bullet\quad&
S\in\mathbf{C}^{\nu}([2h,3h]),\\
\bullet\quad&
S(s)>0,
\qquad s\in(2h,3h),\\
\bullet\quad&
S^{(k)}(2h)=S^{(k)}(3h)=0,
\qquad k=0,1,\ldots,\nu.
\end{aligned}
\label{rou}
\end{equation}
Define a generalized controllability-observability Gramian
\begin{equation}
G = \int_{2h}^{3h} S(s)\mathit{\Gamma}_{\mathrm{d}}(s)
Q_{\mathrm{d}} \mathit{\Gamma}_{\mathrm{d}}^{\mathrm{T}}(s) \mathrm{d}s
\in \mathbb{R}^{n^{2}\times n^{2}},\label{eq:multiple_delay_Gramian}%
\end{equation}
where \begin{equation}
Q_{\mathrm{d}}
=
Q_{\mathrm{d}}^{\mathrm{T}}
\in
\mathbb{R}^{n_{0}mp\times n_{0}mp},
\qquad
Q_{\mathrm{d}}>0
\label{eq:multiple_delay_weighting_matrix}
\end{equation}
is a constant weighting matrix.
\begin{lemma}
\label{thm:multiple_delay_Gramian} Suppose that $(A,B)$ is controllable and
that the selected artificial delays satisfy (\ref{eq:sampled_observability_rank}).
Then the Gramian $G$ defined in (\ref{eq:multiple_delay_Gramian}) is
positive definite, namely, $G>0$.

\end{lemma}

See Appendix A2 for proof.

For $\theta\in[2h,3h)$, define the stacked gain vector
\begin{equation}
\ell_{\mathrm d}(\theta)
=-S(\theta)Q_{\mathrm d}\mathit{\Gamma}_{\mathrm d}^{\mathrm T}(\theta)
G^{-1}\operatorname{vec}(\mathrm e^{Ah})
\in\mathbb R^{n_0mp}.
\label{eq:multiple_delay_gain}
\end{equation}By partitioning $\ell_{\mathrm d}(\theta)$ into $n_0$ blocks as
\[
\ell_{\mathrm d}(\theta)
=\operatorname{col}\bigl(\ell_1(\theta),\ldots,\ell_{n_0}(\theta)\bigr),
\qquad \ell_i(\theta)\in\mathbb R^{mp},
\]
we denote the corresponding gain matrices by
\begin{equation}
L_i(\theta)=\operatorname{vec}^{-1}\!\bigl(\ell_i(\theta)\bigr)
\in\mathbb R^{m\times p},
\qquad i=1,\ldots,n_0.
\label{eq:multiple_delay_gain_matrix}
\end{equation}

With the above preparations, we can establish the following result.

\begin{theorem}
\label{thm:multiple_delay_PT} Consider system (\ref{eq:MIMO_plant}) under
Assumption~\ref{assumption_minimal}. Given $T>0$, choose $h=T/3$ and
$h_i\in[h,2h]$, $i=1,\ldots,n_0$, such that
(\ref{eq:sampled_observability_rank}) holds. Consider the delayed SOF with time-varying gains
\begin{equation}
u(t)=g(t)\sum_{i=1}^{n_0}L_i(\theta)y(t-h_i),
\qquad t\geq 0,
\label{eq:multiple_delay_controller}
\end{equation}
where $\theta$, $g(t)$, and $L_i(\theta)$ are defined in (\ref{eq:cycle_index_phase}), (\ref{eq:act_and_wait_function}), and (\ref{eq:multiple_delay_gain_matrix}), respectively.  Then the closed-loop system (\ref{eq:MIMO_plant}) and (\ref{eq:multiple_delay_controller}) is prescribed-time stable with settling-time bound
$T=3h$. Moreover, $u(t)=0$, $t\geq3h.$
\end{theorem}

\begin{proof}
For $t\in \lbrack0,2h]$, we have $g(t)=0$, and hence the control input
vanishes. Therefore, $x(t)=\mathrm{e}^{At}x(0)$. In particular,
$x(2h)=\mathrm{e}^{2Ah}x(0)$.

For $t\in[2h,3h]$ and $h_{i}\in[h,2h]$, we have $t-h_{i}\in[0,2h]$. Since the
system evolves freely on $[0,2h]$, it follows that
\begin{align*}
x(t-h_{i}) = \mathrm{e}^{A(t-h_{i})}x(0)= \mathrm{e}^{A(t-h_{i}-2h)}x(2h),
\end{align*}
and thus
\[
y(t-h_{i}) = C\mathrm{e}^{A(t-h_{i}-2h)}x(2h).
\]

By the variation-of-constants formula, we can obtain
\begin{align*}
&x(3h)= \\
&  \left[  \mathrm{e}^{Ah} + \sum_{i=1}^{n_{0}} \int_{2h}^{3h}
\mathrm{e}^{A(3h-s)} BL_{i}(s)C \mathrm{e}^{A(s-h_{i}-2h)} \mathrm{d}s
\right] x(2h).
\end{align*}

We next show that the matrix inside the square brackets is zero. Using the
linearity of the vectorization operator and the identity
\[
\operatorname{vec}(XYZ) = (Z^{\mathrm{T}}\otimes X)\operatorname{vec}(Y),
\]
we obtain, for every $i\in \{1,\ldots,n_{0}\}$,
\begin{align*}
& \operatorname{vec}\! \left(  \mathrm{e}^{A(3h-s)} BL_{i}(s)C \mathrm{e}%
^{A(s-h_{i}-2h)} \right) \\
={} &  \left[  \left(  C\mathrm{e}^{A(s-h_{i}-2h)} \right) ^{\mathrm{T}}
\otimes \mathrm{e}^{A(3h-s)}B \right]  \operatorname{vec}\! \left(
L_{i}(s)\right) \\
={} &  \left[  \mathrm{e}^{A^{\mathrm{T}}(s-h_{i}-2h)}C^{\mathrm{T}}
\otimes \mathrm{e}^{A(3h-s)}B \right] \ell_{i}(s)\\
={} &  \mathit{\Gamma}_{\mathrm{d},i}(s)\ell_{i}(s).
\end{align*}
Consequently, by the definitions of $\mathit{\Gamma}_{\mathrm{d}}(s)$ and
$\ell_{\mathrm{d}}(s)$,
\begin{align*}
\mathit{\Gamma}_{\mathrm{d}}(s)\ell_{\mathrm{d}}(s)  & =
\begin{bmatrix}
\mathit{\Gamma}_{\mathrm{d},1}(s) & \cdots & \mathit{\Gamma}_{\mathrm{d}%
,n_{0}}(s)
\end{bmatrix}
\begin{bmatrix}
\ell_{1}(s)\\
\vdots \\
\ell_{n_{0}}(s)
\end{bmatrix}
\\
& = \sum_{i=1}^{n_{0}} \mathit{\Gamma}_{\mathrm{d},i}(s)\ell_{i}(s)\\
& = \sum_{i=1}^{n_{0}} \operatorname{vec}\! \left(  \mathrm{e}%
^{A(3h-s)} BL_{i}(s)C \mathrm{e}^{A(s-h_{i}-2h)} \right) .
\end{align*}
It follows from the linearity of the vectorization operator and integration
that
\begin{align*}
& \operatorname{vec}\! \left(  \sum_{i=1}^{n_{0}} \int_{2h}^{3h}
\mathrm{e}^{A(3h-s)} BL_{i}(s)C \mathrm{e}^{A(s-h_{i}-2h)} \mathrm{d}s
\right) \\
={} &  \int_{2h}^{3h} \sum_{i=1}^{n_{0}} \operatorname{vec}\! \left(
\mathrm{e}^{A(3h-s)} BL_{i}(s)C \mathrm{e}^{A(s-h_{i}-2h)} \right)
\mathrm{d}s\\
={} &  \int_{2h}^{3h} \mathit{\Gamma}_{\mathrm{d}}(s) \ell_{\mathrm{d}}(s)
 \mathrm{d}s.
\end{align*}

Substituting (\ref{eq:multiple_delay_gain}) into the preceding identity gives
\begin{align*}
& \int_{2h}^{3h} \mathit{\Gamma}_{\mathrm{d}}(s) \ell_{\mathrm{d}}(s)
 \mathrm{d}s\\
={} &  -\int_{2h}^{3h} S(s) \mathit{\Gamma}_{\mathrm{d}}(s) Q_{\mathrm{d}}
\mathit{\Gamma}_{\mathrm{d}}^{\mathrm{T}}(s) G^{-1}
\operatorname{vec}\! \left( \mathrm{e}^{Ah}\right) \mathrm{d}s\\
={} &  -\left[  \int_{2h}^{3h} S(s) \mathit{\Gamma}_{\mathrm{d}}(s)
Q_{\mathrm{d}} \mathit{\Gamma}_{\mathrm{d}}^{\mathrm{T}}(s)  \mathrm{d}s
\right]  G^{-1} \operatorname{vec}\! \left( \mathrm{e}^{Ah}\right)
\\
={} &  -GG^{-1} \operatorname{vec}\! \left(
\mathrm{e}^{Ah}\right)
=  -\operatorname{vec}\! \left( \mathrm{e}^{Ah}\right) ,
\end{align*}
where (\ref{eq:multiple_delay_Gramian}) has been used in the third equality.
By the injectivity of the vectorization operator, we obtain
\[
\sum_{i=1}^{n_{0}} \int_{2h}^{3h} \mathrm{e}^{A(3h-s)} BL_{i}(s)C
\mathrm{e}^{A(s-h_{i}-2h)} \mathrm{d}s = -\mathrm{e}^{Ah}.
\]
Consequently, $x(3h) =0$.

We finally show that the state and control input remain zero after $3h$.
During the next waiting interval $[3h,5h]$, the control input is zero. Since
$x(3h)=0$, we obtain $x(t)=0$, $t\in[3h,5h]$. During the subsequent act
interval $[5h,6h]$, for every $t\in[5h,6h]$ and $h_{i}\in[h,2h]$, one has
$t-h_{i}\in[3h,5h]$. Thus, $y(t-h_{i})=0$, $i=1,\ldots,n_{0}$, and
hence $u(t)=0$, $x(t)=0$, $t\in[5h,6h]$. Repeating this argument over all
subsequent wait and act intervals gives $x(t)=0$ and $u(t)=0$ for all $t\geq3h$.

Finally, we prove that there exists a
constant $M>0$ such that
$\sup_{0\leq t\leq3h}\|x(t)\|\leq M\|x(0)\|$.
Since $S(\cdot)$ and
$\mathit{\Gamma}_{\mathrm d}(\cdot)$ are continuous
on $[2h,3h]$ and $G>0$,
the gain $\ell_{\mathrm d}(\cdot)$ is bounded.
It is trivial to prove the above estimate. Together with $x(t)=0$ for all
$t\geq3h$, this establishes Lyapunov stability and
convergence by the prescribed time $T=3h$.
\end{proof}

\section{Feedback with a Current Output and a Single Output Delay}
\label{Ins_output}
In this section, we design a stabilizing feedback controller
that uses the current measurement and a single delayed
measurement, yielding a simpler controller structure (similar to \cite{zhoutac}).
However, the controller gain depends on a matrix $H$
and vectors $v$ and $w$
(see (\ref{eq:MIMO_rank_one_gain}) and
Remark~\ref{remark_1} below), whose existence is guaranteed
under Assumption~\ref{assumption_minimal}
but for which explicit expressions are not provided
in the MIMO case.
In the SISO case, the gain can be readily determined,
as explained later in Remark~\ref{rem:SISO_single_delay}.

Let $K\in\mathbb R^{m\times p}$ be a constant
current-output gain to be designed. Define
\begin{equation}
F_K=A+BKC\in\mathbb R^{n\times n},
\label{eq:Fq}
\end{equation}
and
\begin{equation*}
\begin{aligned}
A_K&=A^{\mathrm T}\otimes I_n-I_n\otimes F_K
\in\mathbb R^{n^2\times n^2},\\
B_C&=C^{\mathrm T}\otimes B
\in\mathbb R^{n^2\times mp}.
\end{aligned}
\end{equation*}
For $s\in[2h,3h]$, define
\begin{align}
\mathit{\Gamma}_K(s)
&:=\mathrm e^{A^{\mathrm T}(s-3h)}C^{\mathrm T}
\otimes\mathrm e^{F_K(3h-s)}B\nonumber\\
&=\mathrm e^{A_K(s-3h)}B_C
\in\mathbb R^{n^2\times mp}.
\label{eq:single_delay_lifted_kernel}
\end{align}
Define a generalized controllability-observability Gramian
\begin{equation}
G_{K}=\int_{2h}^{3h}S(s)\mathit{\Gamma}_K(s)
Q\mathit{\Gamma}_K^{\mathrm T}(s)\mathrm ds,
\label{eq:MIMO_lifted_Gramian}
\end{equation}
where \begin{equation}
Q=Q^{\mathrm T}\in\mathbb R^{mp\times mp},\qquad Q>0,
\label{eq:single_delay_weighting_matrix}
\end{equation}
and $S$ is defined in (\ref{rou}). The controllability
matrix of the pair $(A_K,B_C)$ is
\begin{equation}
\mathcal{K}_{K}=
\begin{bmatrix}
B_C&A_KB_C&\cdots&A_K^{n^2-1}B_C
\end{bmatrix}
\in\mathbb R^{n^2\times n^2mp}.
\label{eq:lifted_controllability_matrix}
\end{equation}

\begin{lemma}
\label{thm:MIMO_constant_Gramian}
Consider system (\ref{eq:MIMO_plant}) under
Assumption~\ref{assumption_minimal}. Then there exists a
constant matrix $K\in\mathbb R^{m\times p}$ such that
$\operatorname{rank}\mathcal{K}_{K}=n^2$.
More precisely, there exist a constant matrix
$H\in\mathbb R^{m\times p}$ and constant vectors
$v\in\mathbb R^m$ and $w\in\mathbb R^p$ such that,
for the family
\begin{equation}
K_\alpha=H+\alpha vw^{\mathrm T},\qquad\alpha\in\mathbb R,
\label{eq:MIMO_rank_one_gain}
\end{equation}
the set of rank-deficient parameter values
\[
\mathcal E_{H,v,w}:=
\left\{\alpha\in\mathbb R:
\operatorname{rank}\mathcal{K}_{K_\alpha}<n^2\right\}
\]
is finite. Moreover, for any constant matrix $K$ satisfying
$\operatorname{rank}\mathcal{K}_{K}=n^2$, the Gramian
defined in (\ref{eq:MIMO_lifted_Gramian}) satisfies $G_{K}>0$.
\end{lemma}

See Appendix A3 for proof.

\begin{remark}
\label{remark_1}
The proof of Lemma~\ref{thm:MIMO_constant_Gramian}
provides a way to select $K$. First, select $H$ so that
$F_H=A+BHC$ has $n$ distinct eigenvalues, which is always
possible under Assumption~\ref{assumption_minimal}
\cite{davison73tac}. Next, choose $v$ and $w$ so that
$(F_H,Bv)$ is controllable and $(F_H,w^{\mathrm T}C)$ is
observable. Finally, choose $\alpha$ outside the finite set
of rank-deficient parameter values and set $K=K_\alpha$. The condition
$\operatorname{rank}\mathcal{K}_{K}=n^2$ is directly
verifiable. This construction does not require a large
gain; it only requires choosing $\alpha$ outside the finite set of rank-deficient parameter values.
\end{remark}

We next present the prescribed-time stabilization result
for system (\ref{eq:MIMO_plant}).
\begin{theorem}
\label{thm:MIMO_exact_PT}
Consider system \eqref{eq:MIMO_plant} under
Assumption~\ref{assumption_minimal}. Given $T>0$, set $h=T/3$
and choose a constant matrix $K\in\mathbb R^{m\times p}$
such that $\operatorname{rank}\mathcal{K}_K=n^2$,
where $\mathcal{K}_K$ is defined in
\eqref{eq:lifted_controllability_matrix}.
By Lemma~\ref{thm:MIMO_constant_Gramian}, such a
matrix always exists, and the Gramian $G_K$ defined in
\eqref{eq:MIMO_lifted_Gramian} is positive definite.
Denote
$L_K(\theta)=\operatorname{vec}^{-1}\!\left(\ell_K(\theta)\right)
\in\mathbb R^{m\times p}$, with
\begin{equation}
\ell_K(\theta)
=-S(\theta)Q\mathit{\Gamma}_K^{\mathrm T}(\theta)
G_K^{-1}\operatorname{vec}\!\left(\mathrm e^{F_Kh}\right)
\in\mathbb R^{mp},
\label{eq:MIMO_lifted_gain}
\end{equation}
where $\theta$, $S(\theta)$, $F_K$, $\mathit{\Gamma}_K(\theta)$, and $Q$ are
defined in \eqref{eq:cycle_index_phase}, \eqref{rou}, \eqref{eq:Fq},
\eqref{eq:single_delay_lifted_kernel}, and \eqref{eq:single_delay_weighting_matrix}, respectively.
Consider the delayed SOF with time-varying gain $L_K(\theta)$
\begin{equation}
u(t)=g(t)Ky(t)+g(t)L_K(\theta)y(t-h),\qquad t\geq0,
\label{eq:MIMO_act_wait_controller}
\end{equation}
where $g(t)$ is defined in \eqref{eq:act_and_wait_function}. Then the closed-loop system \eqref{eq:MIMO_plant} and \eqref{eq:MIMO_act_wait_controller} is prescribed-time stable
with settling-time bound $T=3h$. Moreover, $u(t)=0$ for all $t\geq3h$.
\end{theorem}

\begin{proof}
For $t\in[0,2h]$, we have $g(t)=0$. Therefore,
\begin{equation}
x(t)=\mathrm{e}^{At}x(0),\qquad t\in[0,2h].
\label{eq:MIMO_waiting_state}%
\end{equation}
For $t\in(2h,3h)$, the closed-loop system is
\[
\dot{x}(t)=F_K x(t)
+BL_{K}(t)y(t-h).
\]
Since $t-h\in[h,2h]$, it follows from
(\ref{eq:MIMO_waiting_state}) that
\[
y(t-h)=C\mathrm{e}^{A(t-h)}x(0)
=C\mathrm{e}^{A(t-3h)}x(2h).
\]
By the variation-of-constants formula,
\begin{align*}
x(3h)
={}&\mathrm{e}^{F_K h}x(2h)\\
&+\int_{2h}^{3h}\mathrm{e}^{F_K(3h-s)}
BL_{K}(s)C\mathrm{e}^{A(s-3h)}
\mathrm{d}sx(2h).
\end{align*}
Using the Kronecker--vectorization identity
$
\operatorname{vec}(XYZ)
=
\left(Z^{\mathrm{T}}\otimes X\right)\operatorname{vec}(Y)$,
with $X=\mathrm{e}^{F_K(3h-s)}B$, $Y=L_{K}(s)$, and $Z=C\mathrm{e}^{A(s-3h)}$,
we obtain
\begin{align*}
&\operatorname{vec}\!\left(
\mathrm{e}^{F_K(3h-s)}
BL_{K}(s)C\mathrm{e}^{A(s-3h)}
\right)\\
={}&
\left[
\left(
C\mathrm{e}^{A(s-3h)}
\right)^{\mathrm{T}}
\otimes
\mathrm{e}^{F_K(3h-s)}B
\right]
\operatorname{vec}\!\left(L_{K}(s)\right)\\
={}&
\left[
\mathrm{e}^{A^{\mathrm{T}}(s-3h)}C^{\mathrm{T}}
\otimes
\mathrm{e}^{F_K(3h-s)}B
\right]
\ell_{K}(s)\\
={}&
\mathit{\Gamma}_{K}(s)
\ell_{K}(s).
\end{align*}
Using this identity and (\ref{eq:MIMO_lifted_gain}), we obtain
\begin{align*}
&\operatorname{vec}\!\left(
\int_{2h}^{3h}\mathrm{e}^{F_K(3h-s)}
BL_{K}(s)C\mathrm{e}^{A(s-3h)}\mathrm{d}s
\right)\\
={}&\int_{2h}^{3h}\mathit{\Gamma}_{K}(s)
\ell_{K}(s)\,\mathrm{d}s\\
={}&-\left[
\int_{2h}^{3h}S(s)\mathit{\Gamma}_{K}(s)
Q\mathit{\Gamma}_{K}^{\mathrm{T}}(s)\mathrm{d}s
\right]G_{K}^{-1}
\operatorname{vec}\!\left(\mathrm{e}^{F_K h}\right)\\
={}&-\operatorname{vec}\!\left(
\mathrm{e}^{F_K h}\right).
\end{align*}
Since the vectorization operator is injective, we obtain
\[
\int_{2h}^{3h}\mathrm{e}^{F_K(3h-s)}
BL_{K}(s)C\mathrm{e}^{A(s-3h)}\,\mathrm{d}s
=-\mathrm{e}^{F_K h}.
\]
Consequently, $x(3h)=0$.

Since $x(3h)=0$, the state remains zero throughout the
next waiting interval $[3h,5h]$. During the next act interval,
the delayed state is sampled from $[4h,5h]$ and is therefore
zero. Repeating this argument gives $x(t)=0$ and $u(t)=0$
for all $t\geq3h$.

Finally, $K$ is constant and $L_K(\cdot)$ is bounded on
$[2h,3h]$. The variation-of-constants formula therefore gives
a constant $M>0$, independent of $x(0)$, such that
$\sup_{0\leq t\leq3h}\|x(t)\|\leq M\|x(0)\|$.
Together with $x(t)=0$ for $t\geq3h$, this establishes
Lyapunov stability and convergence by the prescribed time
$T=3h$.
\end{proof}

\begin{remark}
\label{rem:SISO_single_delay}
For SISO systems, the constant current-output gain is
a scalar $K=q\in\mathbb R$, and
$\mathcal K_q\in\mathbb R^{n^2\times n^2}$.
The preceding MIMO result guarantees the existence of
$q_0$ such that $\det\mathcal K_{q_0}\neq0$.
Since the entries of $\mathcal K_q$ are polynomials in $q$,
$\det\mathcal K_q$ is a polynomial that is not identically
zero. Consequently, the set of rank-deficient gain values
\[
\mathcal E
:=
\{q\in\mathbb R:\det\mathcal K_q=0\}
\]
is finite.
By Lemma~\ref{thm:MIMO_constant_Gramian} and
Theorem~\ref{thm:MIMO_exact_PT}, controller
(\ref{eq:MIMO_act_wait_controller}) with $K=q$
achieves prescribed-time stabilization for every
$q\notin\mathcal E$.
In particular, all sufficiently large $|q|$ are admissible.
\end{remark}
\begin{remark}
The controller (\ref{eq:MIMO_act_wait_controller}) can be
modified by replacing $g(t)$ with
\[
g_0(t)=
\begin{cases}
0, & \theta\in[0,h],\\
1, & \theta\in(h,2h),
\end{cases}
\]
and adjusting the cycle length, local phase, and the
corresponding gain formulas accordingly, yielding a
settling time of $2h$, as in \cite{zhoutac}. We retain the design with $g(t)$
for consistency with the previous section.
\end{remark}
\begin{remark}
\label{rem:terminal_cancellation_principle}
The two controller designs in
(\ref{eq:multiple_delay_controller}) and
(\ref{eq:MIMO_act_wait_controller})
share the same principle:
the output-feedback gains are chosen so that the accumulated
contribution of the delayed feedback cancels the free
terminal transition, namely,
$\mathrm e^{Ah}x(2h)$ and
$\mathrm e^{F_Kh}x(2h)$, respectively,
for every state $x(2h)$.

To see this, define
$R_i(s)=C\mathrm e^{A(s-h_i-2h)}$.
During the act interval,
$y(s-h_i)=R_i(s)x(2h)$, and vectorization gives
\[
\operatorname{vec}\!(
\mathit{\Psi}(3h,s)BL_i(s)R_i(s)
)
=
\underbrace{(
R_i^{\mathrm T}(s)\otimes\mathit{\Psi}(3h,s)B
)}_{\text{lifted kernel}}
\operatorname{vec}(L_i(s)),
\]
where $\mathit{\Psi}(3h,s)$ is the state transition
from $s$ to $3h$ during the act interval.
Thus, the columns of the lifted kernel specify the
directions in which the gain can modify the terminal map.

For the multiple-delay design,
$\mathit{\Psi}(3h,s)=\mathrm e^{A(3h-s)}$.
The lifted blocks and their horizontal concatenation are defined in
(\ref{eq:multiple_delay_lifted_block}) and
(\ref{eq:multiple_delay_lifted_matrix}), respectively.
Different delays provide different historical-output
matrices $R_i(s)$ and hence additional groups of
column directions.

For the current-plus-single-delay design, only the
delay $h$ is used, giving $R(s)=C\mathrm e^{A(s-3h)}$.
The corresponding lifted kernel $\mathit{\Gamma}_K(s)$
is defined in (\ref{eq:single_delay_lifted_kernel}).
Here, the constant current-output gain $K$ modifies
the state transition to $\mathrm e^{F_K(3h-s)}$
and thereby changes the column directions of the
lifted kernel.

Under the respective conditions established above,
the columns of $\mathit{\Gamma}_{\mathrm d}(s)$ or
$\mathit{\Gamma}_K(s)$, taken over $s\in(2h,3h)$,
span $\mathbb R^{n^2}$.
The associated Gramian is therefore positive definite,
allowing the delayed-output gains to cancel the free
terminal transition and yield $x(3h)=0$.
\end{remark}

\section{Measurement Delays and Sampled-Data Delayed Implementation}
\label{sec:measurement_delays}

In this section, we extend the multiple-delay design
to systems with known time-varying measurement delays
and unknown constant measurement delays (for nilpotent $A$), and present
a sampled-output implementation.

\subsection{Known Time-Varying Measurement Delays}
\label{subsec:known_tv_measurement_delays}

Consider the LTI system
\begin{equation}
\left\{
\begin{aligned}
\dot{x}(t)&=Ax(t)+Bu(t),\\
y(t)&=
\begin{cases}
Cx(t-\tau(t)), & t-\tau(t)\geq0,\\
0, & t-\tau(t)<0,
\end{cases}
\end{aligned}
\right.
\label{eq:known_tv_delay_system}
\end{equation}
where $A$, $B$, $C$, $x(t)$, $u(t)$, and $y(t)$ have
the same dimensions as in (\ref{eq:MIMO_plant}), and
the measurement delay $\tau(t)$ is known, bounded,
and piecewise continuous, with
$0\leq\tau(t)\leq\tau_M$.

\begin{remark}(Constant Input Delay)
\label{rem:unknown_constant_input_delay}
The formulation (\ref{eq:known_tv_delay_system}) also
accommodates a known constant input delay $\tau_0\geq0$.
Consider the same output equation with
$\dot x(t)=Ax(t)+Bu(t-\tau_0)$ and $u(t)=0$ for $t<0$.
Following Remark~6 of~\cite{zhang26tac}, set
$z(t)=x(t+\tau_0)$. Then
\[
\dot z(t)=Az(t)+Bu(t),\qquad
y(t)=Cz(t-\tau(t)-\tau_0),
\]
where the output identity holds for
$t-\tau(t)\geq\tau_0$.
The effective measurement delay is therefore
$\widetilde\tau(t)=\tau_0+\tau(t)$, with upper bound
$\widetilde\tau_M=\tau_0+\tau_M$.
With the parameter choices below based on
$\widetilde\tau_M$, every delayed measurement used
during an act interval satisfies
$t-h_i-\tau(t-h_i)\geq\tau_0$.
Hence, the initial output convention does not affect
the controller, and the measurement-delay design below
applies using $\widetilde\tau(t)$, with the settling-time
bound increased by $\tau_0$ for the original system.
For nilpotent $A$, the same transformation also handles
an unknown constant input delay via the delay-independent
measurement-delay design in
Section~\ref{subsec:unknown_constant_measurement_delays}.
\end{remark}

For the design developed below, the settling-time bound
can be assigned a priori subject to $T>3\tau_M$ and
therefore cannot be chosen arbitrarily small, consistent
with the delay-dependent restrictions on convergence times
discussed in \cite{zhoutac,zhang26tac,zhang26auto}.
Since the bound is independent of the initial condition,
the resulting stability property is fixed-time stability
in the sense of Definition~\ref{deff}.

The stabilization problem is formulated as follows.

\begin{problem}
Consider system (\ref{eq:known_tv_delay_system}) under
Assumption~\ref{assumption_minimal}.
Given a time bound $T>3\tau_M$, the objective is to design
an artificial-delay-based time-varying SOF
controller such that the closed-loop system is fixed-time
stable with settling-time bound $T$.
\end{problem}

Let $h=T/3>\tau_M$ and choose $n_{0}$ artificial output delays
satisfying
\begin{equation}
h\leq h_1<h_2<\cdots<h_{n_{0}}\leq2h-\tau_M.
\label{eq:known_tv_hi}
\end{equation}
For $s\in[2h,3h]$, define
\[
\mathit{\Gamma}_{\mathrm v,i}(s)
:=
\mathrm e^{A^{\mathrm T}(s-h_i-2h-\tau(s-h_i))}
C^{\mathrm T}
\otimes
\mathrm e^{A(3h-s)}B
\in\mathbb R^{n^2\times mp},
\]
where $i=1,\ldots,n_{0}$, and collect these block matrices as
\begin{equation}
\mathit{\Gamma}_{\mathrm v}(s)
:=
\begin{bmatrix}
\mathit{\Gamma}_{\mathrm v,1}(s)&
\cdots&
\mathit{\Gamma}_{\mathrm v,n_{0}}(s)
\end{bmatrix}
\in\mathbb R^{n^2\times n_{0}mp}.
\label{eq:known_tv_lifted_matrix}
\end{equation}
Define
\begin{equation}
G_{\mathrm v}
=
\int_{2h}^{3h}
S(s)\mathit{\Gamma}_{\mathrm v}(s)
Q_{\mathrm v}
\mathit{\Gamma}_{\mathrm v}^{\mathrm T}(s)\,\mathrm ds,
\label{eq:known_tv_Gramian}
\end{equation}
where
\begin{equation}
Q_{\mathrm v}
=
Q_{\mathrm v}^{\mathrm T}
\in\mathbb R^{n_{0}mp\times n_{0}mp},
\qquad Q_{\mathrm v}>0.
\label{eq:known_tv_weighting_matrix}
\end{equation}
Define the effective delayed observability matrix
\begin{equation*}
\mathcal O_{\boldsymbol h,\tau}(s)
=
\operatorname{col}\!\left(
C\mathrm e^{-A(h_1+\tau(s-h_1))},
\ldots,
C\mathrm e^{-A(h_{n_0}+\tau(s-h_{n_0}))}
\right).
\end{equation*}

\begin{lemma}
\label{thm:known_tv_Gramian_positive}
Suppose that $(A,B)$ is controllable and that there
exists a nonempty open interval
$\mathcal I\subset(2h,3h)$ on which all the functions
$s\mapsto\tau(s-h_i)$, $i=1,\ldots,n_{0}$, are continuous and
\begin{equation}
\operatorname{rank}
\bigl(\mathcal O_{\boldsymbol h,\tau}(s)\bigr)=n,
\qquad s\in\mathcal I.
\label{eq:effective_observability_rank}
\end{equation}
Then the Gramian defined in (\ref{eq:known_tv_Gramian})
is positive definite, namely, $G_{\mathrm v}>0$.
\end{lemma}

See Appendix A4 for proof.

The rank condition in
Lemma~\ref{thm:known_tv_Gramian_positive} depends on
the specific delay function $\tau(t)$ and is directly
verifiable when this function is known.
We next give a more conservative condition that is
independent of the particular delay trajectory.

\begin{lemma}
\label{cor:robust_tv_observability}
Suppose that $(A,B)$ is controllable and that the
artificial delays are selected such that
(\ref{eq:known_tv_hi}) and
(\ref{eq:sampled_observability_rank}) hold.
If $A\neq0$, assume that
\begin{equation}
0\leq\tau_M<
\frac{1}{\|A\|}
\ln\!\left(
1+
\frac{\sigma_{\min}(\mathcal O_{\boldsymbol h})}
{\sqrt{n_{0}}\,\|C\|\mathrm e^{2h\|A\|}}
\right),
\label{eq:robust_tv_observability}
\end{equation}
where $\mathcal O_{\boldsymbol h}$ is defined in
(\ref{eq:sampled_observability_matrix}).
If $A=0$, no additional upper bound on $\tau_M$ is required
for the rank condition.
Then, for every piecewise-continuous delay trajectory
satisfying $0\leq\tau(t)\leq\tau_M$, we have
$\operatorname{rank}
\bigl(\mathcal O_{\boldsymbol h,\tau}(s)\bigr)=n$ for all
$s\in[2h,3h]$.
Consequently, $G_{\mathrm v}>0$.
\end{lemma}

\begin{proof}
First consider $A\neq0$.
For every $i$ and $s\in[2h,3h]$, the $i$th block of
$\mathcal O_{\boldsymbol h,\tau}(s)
-\mathcal O_{\boldsymbol h}$ is
\begin{align*}
C\mathrm e^{-A(h_i+\tau(s-h_i))}
-C\mathrm e^{-Ah_i}=
C\mathrm e^{-Ah_i}
\bigl(\mathrm e^{-A\tau(s-h_i)}-I_n\bigr).
\end{align*}
Since $0\leq\tau(s-h_i)\leq\tau_M$, the series expansion
of the matrix exponential gives
\begin{align*}
&\bigl\|\mathrm e^{-A\tau(s-h_i)}-I_n\bigr\|
\leq
\sum_{k=1}^{\infty}
\frac{\|A^k\|\tau(s-h_i)^k}{k!}\\
&\leq
\sum_{k=1}^{\infty}
\frac{\|A\|^k\tau(s-h_i)^k}{k!}=
\mathrm e^{\|A\|\tau(s-h_i)}-1\leq
\mathrm e^{\|A\|\tau_M}-1.
\end{align*}
Hence,
\begin{align*}
\left\|
C\mathrm e^{-A(h_i+\tau(s-h_i))}
-C\mathrm e^{-Ah_i}
\right\|\leq
\|C\mathrm e^{-Ah_i}\|
\bigl(\mathrm e^{\|A\|\tau_M}-1\bigr).
\end{align*}
By (\ref{eq:known_tv_hi}), $h_i\leq2h-\tau_M\leq2h$, so
\[
\|C\mathrm e^{-Ah_i}\|
\leq\|C\|\mathrm e^{\|A\|h_i}
\leq\|C\|\mathrm e^{2h\|A\|}.
\]
Thus, stacking the block rows yields
\begin{align*}
\bigl\|
\mathcal O_{\boldsymbol h,\tau}(s)
-\mathcal O_{\boldsymbol h}
\bigr\|
&\leq
\left(
\sum_{i=1}^{n_{0}}\|C\mathrm e^{-Ah_i}\|^2
\right)^{1/2}
\bigl(\mathrm e^{\|A\|\tau_M}-1\bigr)\\
&\leq
\sqrt{n_{0}}\,\|C\|\mathrm e^{2h\|A\|}
\bigl(\mathrm e^{\|A\|\tau_M}-1\bigr).
\end{align*}

Condition (\ref{eq:robust_tv_observability}) is equivalent to
\[
\mathrm e^{\|A\|\tau_M}-1
<
\frac{\sigma_{\min}(\mathcal O_{\boldsymbol h})}
{\sqrt{n_{0}}\,\|C\|\mathrm e^{2h\|A\|}}.
\]
Consequently,
$
\bigl\|
\mathcal O_{\boldsymbol h,\tau}(s)
-\mathcal O_{\boldsymbol h}
\bigr\|
<
\sigma_{\min}(\mathcal O_{\boldsymbol h})$
for every $s\in[2h,3h]$.
By the singular-value perturbation inequality \cite{horn13},
\begin{align*}
\sigma_{\min}
\bigl(\mathcal O_{\boldsymbol h,\tau}(s)\bigr)
\geq
\sigma_{\min}(\mathcal O_{\boldsymbol h})-
\bigl\|
\mathcal O_{\boldsymbol h,\tau}(s)
-\mathcal O_{\boldsymbol h}
\bigr\|>0.
\end{align*}
Thus,
$\operatorname{rank}
(\mathcal O_{\boldsymbol h,\tau}(s))=n$
throughout $[2h,3h]$.

If $A=0$, then
$\mathrm e^{-A(h_i+\tau(s-h_i))}
=\mathrm e^{-Ah_i}=I_n$.
Therefore,
$\mathcal O_{\boldsymbol h,\tau}(s)
=\mathcal O_{\boldsymbol h}$
for every $s$, and the rank condition follows directly.

Since $\tau$ is piecewise continuous, there exists a
nonempty open interval $\mathcal I\subset(2h,3h)$
on which all the shifted delay functions
$s\mapsto\tau(s-h_i)$, $i=1,\ldots,n_{0}$, are continuous.
The rank condition holds on this interval.
Hence, Lemma~\ref{thm:known_tv_Gramian_positive}
gives $G_{\mathrm v}>0$.
\end{proof}

For $\theta\in[2h,3h)$, define the stacked gain vector
\begin{equation}
\ell_{\mathrm v}(\theta)
=-S(\theta)Q_{\mathrm v}\mathit{\Gamma}_{\mathrm v}^{\mathrm T}(\theta)
G_{\mathrm v}^{-1}\operatorname{vec}(\mathrm e^{Ah})
\in\mathbb R^{n_0mp},
\label{eq:known_tv_gain}
\end{equation}
where $S(\theta)$, $Q_{\mathrm v}$, $\mathit{\Gamma}_{\mathrm v}(\theta)$,
and $G_{\mathrm v}$ are given in \eqref{rou},
\eqref{eq:known_tv_weighting_matrix},
\eqref{eq:known_tv_lifted_matrix}, and
\eqref{eq:known_tv_Gramian}, respectively.
By partitioning $\ell_{\mathrm v}(\theta)$ into $n_0$ blocks as
\[
\ell_{\mathrm v}(\theta)
=\operatorname{col}\bigl(
\ell_{1,\mathrm v}(\theta),\ldots,\ell_{n_0,\mathrm v}(\theta)
\bigr),
\qquad \ell_{i,\mathrm v}(\theta)\in\mathbb R^{mp},
\]
we denote the corresponding gain matrices by
\begin{equation}
L_i(\theta)=\operatorname{vec}^{-1}\!\bigl(\ell_{i,\mathrm v}(\theta)\bigr)
\in\mathbb R^{m\times p},
\qquad i=1,\ldots,n_0.
\label{eq:known_tv_gain_matrix}
\end{equation}

\begin{theorem}
\label{thm:known_tv_delay}
Consider system (\ref{eq:known_tv_delay_system}) under
Assumption~\ref{assumption_minimal}.
Given $T>3\tau_M$, choose
$h=T/3$ and the artificial delays according to
(\ref{eq:known_tv_hi}).
Suppose that the conditions of
Lemma~\ref{thm:known_tv_Gramian_positive} or the
sufficient conditions of
Lemma~\ref{cor:robust_tv_observability} are satisfied.
Consider the delayed SOF with time-varying gains
\begin{equation}
u(t)=g(t)\sum_{i=1}^{n_{0}}L_i(\theta)y(t-h_i),
\qquad t\geq0,
\label{eq:known_tv_controller}
\end{equation}
where $\theta$, $g(t)$, and $L_i(\theta)$ are defined in
\eqref{eq:cycle_index_phase}, \eqref{eq:act_and_wait_function},
and \eqref{eq:known_tv_gain_matrix}, respectively.
Then the closed-loop system
(\ref{eq:known_tv_delay_system}) and
(\ref{eq:known_tv_controller}) is fixed-time stable
with settling-time bound $T=3h$.
Moreover, $x(t)=u(t)=0$ for all $t\geq3h$.
\end{theorem}

\begin{proof}
Since $g(t)=0$ for $t\in[0,2h]$, the system evolves
freely on this interval:
$x(t)=\mathrm e^{At}x(0)$, $t\in[0,2h]$.
For $s\in[2h,3h]$, condition (\ref{eq:known_tv_hi})
and $0\leq\tau(s-h_i)\leq\tau_M$ imply
$0\leq s-h_i-\tau(s-h_i)\leq2h$.
Consequently, each output used during the first act
interval depends on the state in the preceding waiting
interval:
$y(s-h_i)
=
C\mathrm e^{A(s-h_i-2h-\tau(s-h_i))}x(2h)$.

By the variation-of-constants formula,
\begin{align*}
&x(3h)
={}\mathrm e^{Ah}x(2h)+\\
&\left[
\sum_{i=1}^{n_{0}}
\int_{2h}^{3h}
\mathrm e^{A(3h-s)}BL_i(s)C
\mathrm e^{A(s-h_i-2h-\tau(s-h_i))}
\,\mathrm ds
\right]x(2h).
\end{align*}
The Kronecker--vectorization identity gives
\begin{align*}
&\operatorname{vec}\!\left(
\sum_{i=1}^{n_{0}}
\int_{2h}^{3h}
\mathrm e^{A(3h-s)}BL_i(s)C
\mathrm e^{A(s-h_i-2h-\tau(s-h_i))}
\,\mathrm ds
\right)\\
&\qquad=
\int_{2h}^{3h}
\mathit{\Gamma}_{\mathrm v}(s)
\ell_{\mathrm v}(s)\,\mathrm ds.
\end{align*}
Substituting (\ref{eq:known_tv_gain}) and using
(\ref{eq:known_tv_Gramian}) yields
$
\int_{2h}^{3h}
\mathit{\Gamma}_{\mathrm v}(s)
\ell_{\mathrm v}(s)\,\mathrm ds
=
-\operatorname{vec}(\mathrm e^{Ah})$.
Hence, $x(3h)=0$.

During the next waiting interval $[3h,5h]$, the state
remains zero.
For $t\in(5h,6h)$, the delay bounds imply
$3h\leq t-h_i-\tau(t-h_i)\leq5h$,
$i=1,\ldots,n_{0}$.
Thus, all delayed outputs used by the controller vanish.
Repeating this argument yields
$x(t)=u(t)=0$ for all $t\geq3h$.

Since the gains are bounded on $[2h,3h]$, there exists
$M>0$ such that
$
\sup_{0\leq t\leq3h}\|x(t)\|
\leq M\|x(0)\|$.
Therefore, the closed-loop system is fixed-time stable
with settling-time bound $T=3h$.
\end{proof}
\subsection{Sampled-Data Delayed Implementation}
In the context of sampled-data stabilization
\cite{Mazenc26fixed,karafyllis11tac}, the following
corollary extends the measurement-delay design to
sampled-output implementation with time-varying gains.
\begin{corollary}
\label{rem:sampling_delay_rank_preservation}
Consider system (\ref{eq:MIMO_plant}) under
Assumption~\ref{assumption_minimal}, with output samples
\begin{equation}
y(t_k)=Cx(t_k),\qquad
t_k=kd_s,\qquad k\in\mathbb Z_0^+,
\label{eq:sampled_output_measurements}
\end{equation}
where $d_s>0$ is the sampling period.
Let $T=3h$ with $h>d_s$, and choose values
$h_i=\nu_i d_s$, $\nu_i\in\mathbb N$, $i=1,\ldots,n_{0}$,
satisfying
$h\leq h_1<\cdots<h_{n_{0}}\leq2h-d_s$
and (\ref{eq:sampled_observability_rank}).
Consider the time-varying sampled-data SOF with delays $\nu_i$
\begin{equation}
u(t)=g(t)\sum_{i=1}^{n_{0}}L_i(\theta)y(t_{k-\nu_i}),
\qquad t\in[t_k,t_{k+1}),
\label{eq:sampled_output_controller}
\end{equation}
where $g(t)$ is defined in
(\ref{eq:act_and_wait_function}),
$\theta=t-3h\lfloor t/(3h)\rfloor$, and
$L_i(s)=\operatorname{vec}^{-1}(\ell_{i,\mathrm v}(s))
\in\mathbb R^{m\times p}$, $i=1,\ldots,n_{0}$,
with the stacked gain vector constructed according to
(\ref{eq:known_tv_gain}) using
$\tau(t)=t-t_k$, $t\in[t_k,t_{k+1})$, and $\tau_M=d_s$.
Then the closed-loop system
(\ref{eq:MIMO_plant}) and
(\ref{eq:sampled_output_controller}) is fixed-time stable
with settling-time bound $T=3h$.
Moreover, $u(t)=0$ for all $t\geq3h$.
\end{corollary}

\begin{proof}
Define the held output by
$\bar y(t)=y(t_k)$ for $t\in[t_k,t_{k+1})$. Since $\tau(t)=t-t_k\in[0,d_s)$ and $h_i=\nu_i d_s$,
we have $\bar y(t-h_i)=y(t_{k-\nu_i})$ for
$t\in[t_k,t_{k+1})$ whenever the controller is active.
Thus, controller (\ref{eq:known_tv_controller}) takes
the form (\ref{eq:sampled_output_controller}).
The delay $\tau$ is $d_s$-periodic, so
$\tau(s-h_i)=\tau(s)$ whenever the arguments are nonnegative.
Consequently,
$\mathcal O_{\boldsymbol h,\tau}(s)
=
\mathcal O_{\boldsymbol h}\mathrm e^{-A\tau(s)}$.
Since the matrix exponential is nonsingular,
(\ref{eq:sampled_observability_rank}) gives
$
\operatorname{rank}
\bigl(\mathcal O_{\boldsymbol h,\tau}(s)\bigr)
=
\operatorname{rank}
\bigl(\mathcal O_{\boldsymbol h}\bigr)=n$.
Hence, Lemma~\ref{thm:known_tv_Gramian_positive}
gives $G_{\mathrm v}>0$.
Since $h>d_s=\tau_M$ and
$h_i\in[h,2h-\tau_M]$,
Theorem~\ref{thm:known_tv_delay} yields the stated
fixed-time stability and $u(t)=0$ for $t\geq3h$.
\end{proof}

\subsection{Unknown Constant Delays for Nilpotent $A$}
\label{subsec:unknown_constant_measurement_delays}
We consider unknown constant measurement delays.
By the transformation in
Remark~\ref{rem:unknown_constant_input_delay},
the result below also applies to unknown constant
input delays, with an additional $\tau_M$ in the
settling-time bound.
Consider the LTI system
\begin{equation}
\left\{
\begin{aligned}
\dot{x}(t)&=Ax(t)+Bu(t),\\
y(t)&=
\begin{cases}
Cx(t-\tau_{1}), & t-\tau_{1}\geq0,\\
0, & t-\tau_{1}<0,
\end{cases}
\end{aligned}
\right.
\label{eq:unknown_constant_measurement_delay_system}
\end{equation}
where $\tau_{1}\in[0,\tau_M]$ is an unknown constant,
$\tau_M\geq0$ is a known upper bound,
$(A,B)$ is controllable, $(A,C)$ is observable,
and $A^r=0$ for some positive integer $r\leq n$.

Related fixed-time stabilization problems were studied
using state feedback in \cite{zhang26tac} and SOF in \cite{zhang26auto}, with the latter
restricted to systems without internal dynamics
\cite{Zhou23automatica}.
The following corollary extends the SOF result to systems with internal dynamics,
including nonminimum-phase systems.

\begin{corollary}
\label{prop:unknown_constant_measurement_delay}
Consider system
(\ref{eq:unknown_constant_measurement_delay_system}).
Choose $h>\tau_M$ and $n_{0}=r\leq n$ distinct artificial delays
$h_i\in[h,2h-\tau_M]$, $i=1,\ldots,r$.
Then controller (\ref{eq:multiple_delay_controller}),
with the gains defined in (\ref{eq:multiple_delay_gain}),
renders the closed-loop system fixed-time stable
with settling-time bound $3rh$.
Moreover, $x(t)=u(t)=0$ for all $t\geq3rh$.
\end{corollary}

\begin{proof}
Since $A^r=0$, the delayed observability matrix satisfies
\[
\mathcal O_{\boldsymbol h}
=
(V_{\boldsymbol h}\otimes I_p)
\operatorname{col}(C,CA,\ldots,CA^{r-1}),
\]
where
\[
V_{\boldsymbol h}
=
\begin{bmatrix}
1&-h_1&\cdots&\dfrac{(-h_1)^{r-1}}{(r-1)!}\\
\vdots&\vdots&&\vdots\\
1&-h_r&\cdots&\dfrac{(-h_r)^{r-1}}{(r-1)!}
\end{bmatrix}.
\]
The matrix $V_{\boldsymbol h}$ is nonsingular because
the delays are distinct.
Moreover, observability of $(A,C)$ and $A^r=0$ imply
$\operatorname{rank}\bigl(\operatorname{col}(C,CA,\ldots,CA^{r-1})\bigr)=n$.
Thus, $\operatorname{rank}\mathcal O_{\boldsymbol h}=n$,
and the gains in (\ref{eq:multiple_delay_gain})
are well defined.

During each act interval, every measurement used
by the controller refers to the preceding waiting
interval, since
$0\leq s-h_i-\tau_{1}\leq2h$,
$s\in[2h,3h]$.
The terminal-map cancellation in
Theorem~\ref{thm:multiple_delay_PT} therefore gives
\[
x\bigl(3h(l+1)\bigr)
=
\mathit{\Delta}_{\tau_{1}}x(3hl),
\mathit{\Delta}_{\tau_{1}}
=
\mathrm e^{3Ah}(I-\mathrm e^{-A\tau_{1}}),
l\in\mathbb Z_0^+.
\]
Since $A^r=0$, there exists a polynomial
$\mathcal P_{\tau_{1}}$ such that
$I-\mathrm e^{-A\tau_{1}}=A\mathcal P_{\tau_{1}}(A)$.
All factors commute, and hence
$
\mathit{\Delta}_{\tau_{1}}^r
=
\mathrm e^{3rhA}A^r\mathcal P_{\tau_{1}}(A)^r=0$.
Consequently, $x(3rh)=0$.
The subsequent waiting and act intervals preserve
the zero state and yield $u(t)=0$ for $t\geq3rh$.

The bounded gains and the variation-of-constants formula
give a constant $M>0$, independent of $x(0)$, such that
$\sup_{0\leq t\leq3rh}\|x(t)\|
\leq M\|x(0)\|$.
Together with $x(t)=0$ for $t\geq3rh$, this establishes
Lyapunov stability and fixed-time convergence.
\end{proof}

\begin{remark}
\label{rem:Gramian_conditioning}
The lifted directions described in
Remark~\ref{rem:terminal_cancellation_principle}
also determine the numerical conditioning of the
associated Gramians.
Although the conditions established above guarantee
$G>0$ and $G_K>0$, nearly dependent lifted directions
or large differences in their magnitudes may cause
ill-conditioning.
For example, let
$d_1=[1,0]^{\mathrm T}$,
$d_2=[\cos\vartheta,\sin\vartheta]^{\mathrm T}$,
and $D=[d_1,d_2]$, where $0<\vartheta\leq\pi/2$.
Then
\[
\kappa_2(DD^{\mathrm T})
=\frac{\lambda_{\max}(DD^{\mathrm T})}
{\lambda_{\min}(DD^{\mathrm T})}
=\frac{1+\cos\vartheta}{1-\cos\vartheta}
\longrightarrow\infty
\quad\text{as }\vartheta\to0.
\]
Alternatively, taking
$d_1=[1,0]^{\mathrm T}$ and
$d_2=[0,\delta]^{\mathrm T}$, with $0<\delta\leq1$,
gives orthogonal directions with unequal magnitudes,
for which $DD^{\mathrm T}=\operatorname{diag}(1,\delta^2)$
and $\kappa_2(DD^{\mathrm T})=\delta^{-2}$.
The weighting function and weighting matrices also
affect conditioning by changing the relative
contributions of different times and channels.
These considerations motivate the numerical implementation
below, which determines the gain coefficients directly
from terminal-map cancellation while avoiding explicit
Gramian inversion.
\end{remark}
\section{Numerical Implementation of Gramian-Based Gains}
\label{subsec:finite_dimensional_implementation}

The gain functions in Theorem~\ref{thm:multiple_delay_PT}
are explicitly constructed using the Gramian $G$,
as in the Gramian-based constructions in
\cite{Mazenc26fixed,zhoutac}.
As discussed in Remark~\ref{rem:Gramian_conditioning},
$G$ may be ill-conditioned, making the gain computation
sensitive to numerical errors; the resulting feedback
gains may also be large.
Under an additional assumption
(see (\ref{eq:finite_rank_condition}) below),
we therefore present an alternative implementation
that parameterizes the gain functions in a chosen
function space and determines their coefficients
directly from the terminal-map cancellation condition.
This implementation preserves the act-and-wait structure
without explicitly computing $G^{-1}$.

Let $\{\phi_0,\phi_1,\ldots,\phi_{N_b-1}\}$
be a finite number of some basis functions in $\mathbf C^\nu([2h,3h])$ (e.g. functions given by (\ref{eq:shifted_Legendre_basis})). For each artificial delay $h_i$, $i=1,\ldots,n_{0}$,
assume that the corresponding vectorized gain in (\ref{eq:multiple_delay_controller}) can be
represented as
\begin{equation}
\operatorname{vec}\!\left(L_i(\theta)\right)
=
S(\theta)
\sum_{k=0}^{N_b-1}
\phi_k(\theta)a_{i,k},
\qquad \theta \in[2h,3h),
\label{eq:finite_gain_parameterization}
\end{equation}
where $S(\theta)$ is defined by (\ref{rou}) and $a_{i,k}\in\mathbb R^{mp}$ are constant coefficient
vectors to be determined. In particular, the endpoint conditions
imposed on \(S(\theta)\) ensure that the periodically extended gain functions
satisfy the required boundary conditions.

Following the derivation in the proof of
Theorem~\ref{thm:multiple_delay_PT}, define the
terminal map $F$ by
\[
x(3h)=F x(2h).
\]
By the variation-of-constants formula, it is given by
\[
F
=
\mathrm e^{Ah}
+
\sum_{i=1}^{n_{0}}
\int_{2h}^{3h}
\mathrm e^{A(3h-s)}
B L_i(s) C
\mathrm e^{A(s-h_i-2h)}
\mathrm ds.
\]
Using the Kronecker--vectorization identity and the definition of
\(\mathit{\Gamma}_{\mathrm d,i}(s)\), we obtain
\begin{align*}
\operatorname{vec}(F)
={}&
\operatorname{vec}\!\left(\mathrm e^{Ah}\right)
+
\sum_{i=1}^{n_{0}}
\int_{2h}^{3h}
\mathit{\Gamma}_{\mathrm d,i}(s)
\operatorname{vec}\!\left(L_i(s)\right)
\mathrm ds.
\end{align*}
Substituting (\ref{eq:finite_gain_parameterization}) gives
\begin{align*}
\operatorname{vec}(F)
={}&
\operatorname{vec}\!\left(\mathrm e^{Ah}\right)\\
&+
\sum_{i=1}^{n_{0}}
\sum_{k=0}^{N_b-1}
\left[
\int_{2h}^{3h}
S(s)\phi_k(s)
\mathit{\Gamma}_{\mathrm d,i}(s)
\mathrm ds
\right]a_{i,k}.
\end{align*}

For each \(i=1,\ldots,n_{0}\) and \(k=0,\ldots,N_b-1\), define
\begin{equation*}
K_{i,k}
:=
\int_{2h}^{3h}
S(s)\phi_k(s)
\mathit{\Gamma}_{\mathrm d,i}(s)
\mathrm ds
\in\mathbb R^{n^2\times mp}.
\end{equation*}
Collect all coefficient vectors in the order
\[
a
:=
\operatorname{col}
\{
a_{1,0},\ldots,a_{1,N_b-1},
\ldots,
a_{n_{0},0},\ldots,a_{n_{0},N_b-1}
\}
,
\]
with $a\in\mathbb R^{n_{0}mpN_b}$, and define
\begin{equation}
K_b
:=
\begin{bmatrix}
K_{1,0} & \cdots & K_{1,N_b-1} &
\cdots &
K_{n_{0},0} & \cdots & K_{n_{0},N_b-1}
\end{bmatrix}
,
\label{eq:finite_K_matrix}
\end{equation}
with $K_b\in\mathbb R^{n^2\times n_{0}mpN_b}$. The preceding double sum can then be written compactly as \(K_ba\).
Consequently, the terminal map satisfies
\begin{equation}
\operatorname{vec}(F)
=
\operatorname{vec}\!\left(\mathrm e^{Ah}\right)
+
K_ba.
\label{eq:finite_terminal_map}
\end{equation}

Thus, exact cancellation of the one-cycle terminal map is equivalent to the
finite-dimensional linear equation
\begin{equation}
K_ba
=
-\operatorname{vec}\!\left(\mathrm e^{Ah}\right).
\label{eq:finite_terminal_cancellation}
\end{equation}
Suppose that
\begin{equation}
\operatorname{rank}(K_b)=n^2.
\label{eq:finite_rank_condition}
\end{equation}
Since $K_b$ has $n^2$ rows, this condition means
that $K_b$ has full row rank. Consequently,
$K_bK_b^{\mathrm T}$ is positive definite, and
(\ref{eq:finite_terminal_cancellation}) is solvable
for every right-hand side.
If $n_0mpN_b>n^2$, the equation admits infinitely
many solutions, providing additional freedom to
optimize transient performance.
The unique minimum-Euclidean-norm solution is
\begin{equation}
\begin{aligned}
a=-K_b^{\mathrm T}
\left(K_bK_b^{\mathrm T}\right)^{-1}
\operatorname{vec}\!\left(\mathrm e^{Ah}\right)=-K_b^\dagger
\operatorname{vec}\!\left(\mathrm e^{Ah}\right),
\end{aligned}
\label{eq:finite_pseudoinverse_solution}
\end{equation}
where
$K_b^\dagger
=K_b^{\mathrm T}(K_bK_b^{\mathrm T})^{-1}$
is the Moore--Penrose inverse of $K_b$.

Substitution of (\ref{eq:finite_pseudoinverse_solution}) into
(\ref{eq:finite_terminal_map}) gives $\operatorname{vec}(F)=0$
because the full-row-rank condition implies
$K_bK_b^\dagger=I_{n^2}$.
Therefore,
$x(3h)=0$ for every initial condition, which together with the proof of Theorem~\ref{thm:multiple_delay_PT} yields
$x(t)=0$ and $u(t)=0$, $t\geq3h$. In numerical computations, the pseudoinverse is
evaluated using the singular value decomposition
rather than explicitly forming and inverting
$K_bK_b^{\mathrm T}$.

The minimum-norm solution above ensures terminal cancellation
but does not directly optimize transient performance.
To retain a precomputed coefficient vector $a_0$ as closely
as possible while enforcing terminal cancellation, we apply
the following minimum-norm correction.
An initial coefficient vector $a_0$ may be obtained
from a preliminary numerical optimization. Its terminal-map residual is
$r_0
:=
\operatorname{vec}\!\left(\mathrm e^{Ah}\right)+K_ba_0$.
To remove this residual, consider the corrected coefficient vector
$a
=
a_0
-
K_b^\dagger r_0$,
or, equivalently,
\begin{equation}
a
=
a_0
+
K_b^\dagger
\left[
-\operatorname{vec}\!\left(\mathrm e^{Ah}\right)
-K_ba_0
\right].
\label{eq:finite_coefficient_refinement}
\end{equation}
Indeed, using \(K_bK_b^\dagger=I_{n^2}\), we obtain
\begin{align*}
\operatorname{vec}\!\left(\mathrm e^{Ah}\right)+K_ba
&=
r_0-K_bK_b^\dagger r_0=
0.
\end{align*}
Therefore, the correction in
(\ref{eq:finite_coefficient_refinement}) exactly removes the residual of the
discretized terminal-map equation. Moreover, \(-K_b^\dagger r_0\) is the
minimum-Euclidean-norm correction among all coefficient corrections that
eliminate this residual.

For the numerical examples, the shifted Legendre basis is selected as
\begin{equation}
\phi_k(\theta)
=
P_k\!\left(
\frac{2(\theta-2h)}{h}-1
\right),
k=0,1,\ldots,N_b-1,
\label{eq:shifted_Legendre_basis}
\end{equation}
where \(P_k\) denotes the Legendre polynomial of degree \(k\). The affine
mapping
$\theta\longmapsto 2(\theta-2h)/h-1$
maps the act interval \([2h,3h]\) onto \([-1,1]\), on which the Legendre
polynomials are mutually orthogonal. This choice generally provides a better-conditioned representation than monomial basis functions. Nevertheless, the
Legendre basis is only a numerical choice and is not essential to the control
design. Any collection of basis functions in $\mathbf C^\nu([2h,3h])$ may be used, provided that
the corresponding matrix \(K_b\) satisfies
(\ref{eq:finite_rank_condition}).

The condition \(G>0\) in
Lemma~\ref{thm:multiple_delay_Gramian} guarantees the existence of gain
functions in the underlying infinite-dimensional function space, but it does
not imply that an arbitrarily prescribed finite collection of basis functions
is sufficiently rich to satisfy (\ref{eq:finite_rank_condition}).
Consequently, condition (\ref{eq:finite_rank_condition}) must be verified for
the selected basis functions and basis order. Increasing \(N_b\) enlarges the
admissible finite-dimensional function space and may restore the required rank
or reduce the continuous terminal-map approximation error. However, an
excessively large \(N_b\) may also increase the condition number of \(K_b\)
and amplify quadrature and roundoff errors. The basis order, quadrature
accuracy, numerical rank tolerance, and pseudoinverse regularization should
therefore be selected jointly. We illustrate this implementation in
Section~\ref{ex:cart} below,
where the Gramian is ill-conditioned.

\section{Numerical Examples}
\label{Num_ex}
We consider four examples of nonminimum-phase systems from the literature, to which the method of \cite{zhang26auto} cannot be applied. Throughout the numerical examples, following
\cite{zhang26tac,zhang26auto}, we use
$S(t)=\sin^2(\pi(t-2h)/h)$, $t\in[2h,3h]$,
which satisfies (\ref{rou}) with $\nu=1$.
\subsection{SISO: Nonminimum-Phase Systems}

Consider the SISO system \cite{Michiels04ima,zhang22tac}
\begin{equation}
A=
\begin{bmatrix}
0 & 1\\
1 & -0.5
\end{bmatrix},
\qquad
B=
\begin{bmatrix}
0\\
1
\end{bmatrix},
\qquad
C=
\begin{bmatrix}
1 & -1
\end{bmatrix}.
\label{eq:example2_system}
\end{equation}
Its transfer function is
\[
G(s)=C(sI-A)^{-1}B
=\frac{1-s}{s^{2}+\frac12 s-1}.
\]
The system has a right-half-plane transmission zero at $s_{\mathrm z}=1$
and is therefore nonminimum phase.

We consider $T=3\,\mathrm{s}$ and set $h=T/3=1\,\mathrm{s}$.
For the multiple-delay controller \eqref{eq:multiple_delay_controller}
in Theorem~\ref{thm:multiple_delay_PT}, we use two artificial delays
$h_1=1.0215\,\mathrm{s}$ and $h_2=1.9432\,\mathrm{s}$, and set
\[
Q_{\mathrm d}=
\begin{bmatrix}
 0.4472 & -0.7776\\
-0.7776 &  1.5528
\end{bmatrix}.
\]
For the current-plus-single-delay controller
\eqref{eq:MIMO_act_wait_controller}
in Theorem~\ref{thm:MIMO_exact_PT}, we set $K=q$ as in
Remark~\ref{rem:SISO_single_delay}, with $q=-3.3457$.

The controller parameters $h_1$, $h_2$, $Q_{\mathrm d}$, and $q$
have been found by using the Nelder--Mead simplex method, where the
objective combines the worst-case values of the peak state norm,
peak control magnitude, and control energy over all unit-norm
initial states.

For simulation, take $x(0)=[1,-1]^{\mathrm T}$.
The simulations are run over $0\le t\le 20\,\mathrm{s}$ with measured
output $y(t)=Cx(t)+w(t)$. Both proposed controllers are tested with
$w(t)=0$ and with sinusoidal measurement noise
$w(t)=10^{-3}\sin(12t)$, using the same controller parameters in both
cases.

For comparison, in the noise-free case we also simulate
the square-wave dither controller proposed in \cite{zhang22tac},
\[
u(t)=\frac{k}{\varepsilon}
\operatorname{sgn}\!\left(\cos\!\left(\frac{2\pi t}{\varepsilon}\right)\right)
y(t).
\]
We use $k=9$ and $\varepsilon=3.8$, as in \cite{zhang22tac} for this system.

\begin{figure}[t]
\centering
\includegraphics[scale=0.54]{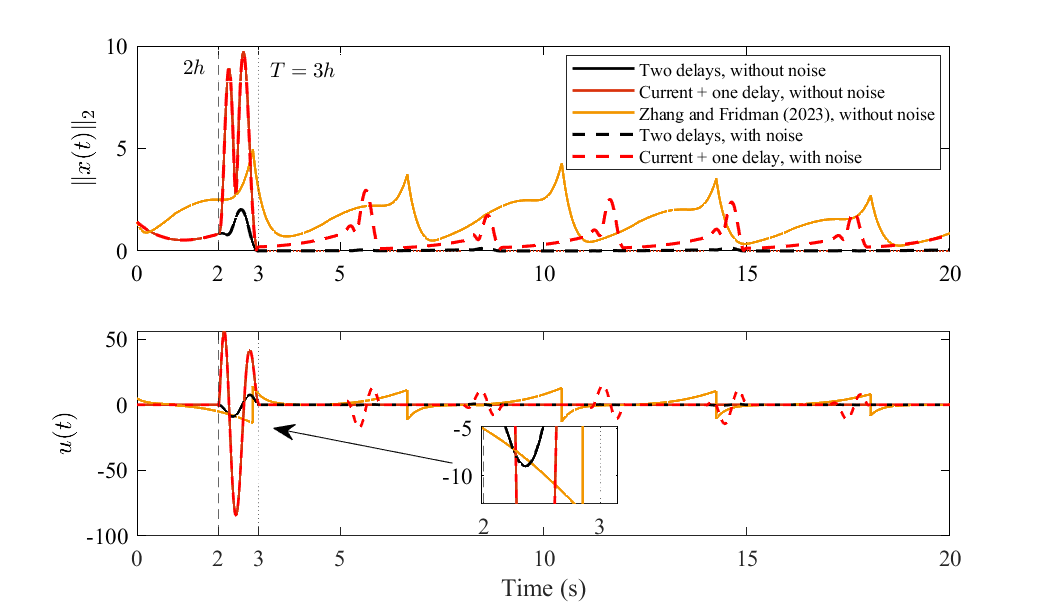}
\caption{State norm and control input for system
\eqref{eq:example2_system}
with and without measurement noise.}
\label{figure1}
\end{figure}

Fig.~\ref{figure1} shows the state norm and the control input.
In the noise-free case, both proposed controllers drive the state to
the origin at the prescribed time $T=3\,\mathrm{s}$. The multiple-delay
controller yields a smaller peak state norm than the
current-plus-single-delay controller and a smaller peak control
magnitude than both the current-plus-single-delay and square-wave
dither controllers. With measurement noise, the multiple-delay
controller exhibits smaller excursions in the state norm than the
current-plus-single-delay controller.

\subsection{MIMO: A Nonminimum-Phase System with a Nilpotent State Matrix}
\label{subsec:MIMO_nilpotent_example}

Consider system (\ref{eq:MIMO_plant}), where the matrices
$A$ and $B$ are taken from \cite{Zhou23automatica},
and the output matrix $C$ is chosen to obtain a
nonminimum-phase system:
\[
A=
\begin{bmatrix}
0&1&-1&2&-1&-1\\
0&-1&0&-1&0&1\\
0&1&0&2&0&-1\\
0&0&0&0&0&0\\
0&0&0&0&0&0\\
0&-1&0&-1&0&1
\end{bmatrix},
B=
\begin{bmatrix}
1&0&0\\
0&-1&0\\
1&1&0\\
0&1&1\\
-1&0&1\\
0&1&1
\end{bmatrix},
\]
and
\[
C=
\begin{bmatrix}
-1&3&-1&7&-4&-3\\
0&2&0&3&-2&-1\\
0&0&0&1&0&-1
\end{bmatrix}.
\]
Clearly, the system is controllable and observable, and has a right-half-plane transmission zero
at $s=1$ and one-dimensional unstable zero dynamics.
Moreover, $A^2\neq0$ and $A^3=0$.

For both controller designs below, we consider $T=9\,\mathrm{s}$
and $T=18\,\mathrm{s}$. For simulation, we take the initial condition
$x(0)=[1,1,1,1,1,1]^{\mathrm T}/\sqrt6$, so that $\|x(0)\|=1$.
\subsubsection{Current-Plus-Single-Delay Output Feedback}

We first consider the system without measurement delays
and apply controller (\ref{eq:MIMO_act_wait_controller})
with $h=T/3$.
To ensure that $F_H$ has six distinct eigenvalues, we select
$H=\operatorname{diag}(1,2,3)$.
Choose $v=[1,0,0]^{\mathrm T}$ and $w=[0,0,1]^{\mathrm T}$,
which gives
$b=Bv=[1,0,1,0,-1,0]^{\mathrm T}$ and
$c=w^{\mathrm T}C=[0,0,0,1,0,-1]$.
Direct calculation yields
$\det[b,F_Hb,\ldots,F_H^5b]=-27648\neq0$ and
$\det[c^{\mathrm T},F_H^{\mathrm T}c^{\mathrm T},
\ldots,(F_H^{\mathrm T})^5c^{\mathrm T}]=256\neq0$.
Hence, $(F_H,b,c)$ is a minimal SISO realization.

We then select $\alpha=1$ so that $F_K=A+BKC$
has six distinct eigenvalues, giving
\[
K=H+\alpha vw^{\mathrm T}
=
\begin{bmatrix}
1&0&1\\
0&2&0\\
0&0&3
\end{bmatrix}.
\]
Since the only distinct eigenvalue of $A$ is zero,
the eigenvalue-difference condition in the proof of
Lemma~\ref{thm:MIMO_constant_Gramian} is satisfied.
Consequently, the lifted pair $(A_K,B_C)$ is controllable
and $G_K>0$.

We take $Q=I_9$ and construct the delayed-output gain
according to (\ref{eq:MIMO_lifted_gain}).
For $T=9\,\mathrm{s}$ and $T=18\,\mathrm{s}$, we select $h=3$ and $h=6$,
respectively.
By Theorem~\ref{thm:MIMO_exact_PT}, the resulting
controllers achieve prescribed-time stabilization
with the corresponding settling-time bounds (see Fig.~\ref{new_figure2}).

\subsubsection{Unknown Constant Measurement Delays}

We next consider system
(\ref{eq:unknown_constant_measurement_delay_system})
with the same matrices $A$, $B$, and $C$,
where the constant measurement delay
$\tau_{1}\in[0,\tau_M]$ is unknown.
Since $A^3=0$, Corollary~\ref{prop:unknown_constant_measurement_delay}
guarantees fixed-time stabilization with a settling-time bound of
$T=9h$ using three artificial delays, collected in the vector
$\boldsymbol{h}=\begin{bmatrix}h_1 & h_2 & h_3\end{bmatrix}^{\mathrm T}$.

For $T=9\,\mathrm{s}$, we select
$h=1$, $\tau_{1}=\tau_M=0.5$,
$\boldsymbol{h}
=\begin{bmatrix}1 & 1.25 & 1.5\end{bmatrix}^{\mathrm T}$.
For $T=18\,\mathrm{s}$, we select
$h=2$, $\tau_{1}=\tau_M=1$, $
\boldsymbol h=
\begin{bmatrix}
2&2.5&3
\end{bmatrix}^{\mathrm T}$.
Both choices satisfy $h>\tau_M$ and
$h\leq h_1<h_2<h_3\leq2h-\tau_M$.

Since $A^3=0$ and the three artificial delays are distinct,
the delayed observability condition follows from the proof of
Corollary~\ref{prop:unknown_constant_measurement_delay}.
We select $Q_{\mathrm d}=I_{27}$ and apply controller
(\ref{eq:multiple_delay_controller}), with the gains
constructed from (\ref{eq:multiple_delay_gain}).
The controller periods are $3h=3$ and $3h=6$,
respectively.
Since the period-to-period state map satisfies
$\mathit{\Delta}_{\tau_{1}}^3=0$, fixed-time convergence is guaranteed
within three periods, giving the corresponding
settling-time bounds $T=9\,\mathrm{s}$ and $T=18\,\mathrm{s}$.

\begin{figure}[t]
\centering
\includegraphics[scale=0.55]{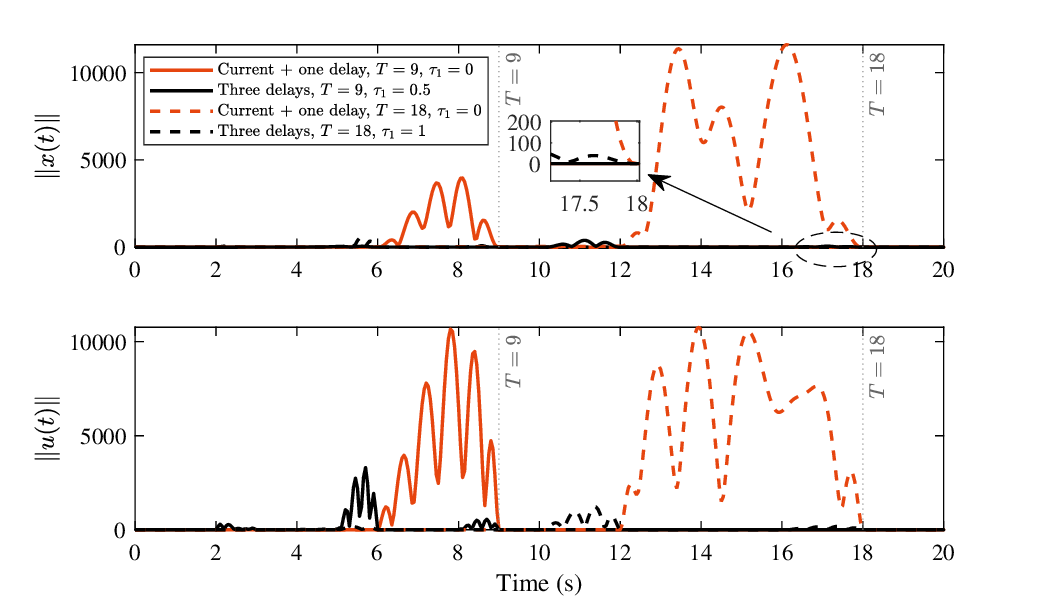}
\caption{State and control-input norms of the nonminimum-phase
MIMO system under the two controller designs with
$T=9\,\mathrm{s}$ and $T=18\,\mathrm{s}$.}
\label{new_figure2}
\end{figure}

The simulation results are shown in Fig.~\ref{new_figure2}.
Both controller designs drive the state to the origin
by the corresponding time bounds $T=9\,\mathrm{s}$ and $T=18\,\mathrm{s}$, and the control inputs
vanish thereafter.
For the selected parameters, the three-delay controller
produces smaller transient state and control peaks
than the current-plus-single-delay controller,
despite the unknown measurement delay.
Although the current-plus-single-delay controller has a
simpler structure, it may produce larger control inputs
and greater transient state peaks.

\subsection{Flexible-Structure System}

To illustrate the applicability of the proposed method
to a physical system, consider the flexible-structure system studied in \cite{qiu93auto}.
We adopt the state-space model reported therein, which is
described by (\ref{eq:MIMO_plant}) with
\[
A=
\begin{bmatrix}
-0.996 & -463 & -97.8 & -12131 & -8.11 & 0\\
1 & 0 & 0 & 0 & 0 & 0\\
0 & 1 & 0 & 0 & 0 & 0\\
0 & 0 & 1 & 0 & 0 & 0\\
0 & 0 & 0 & 1 & 0 & 0\\
0 & 0 & 0 & 0 & 1 & 0
\end{bmatrix},
\]
and
\[
B=[1,0,0,0,0,0]^{\mathrm T},
C=[0,1.65,-0.331,-576,90.6,19080].
\]
Thus, the system order is $n=6$.
The transfer function possesses two right-half-plane zeros,
approximately located at $6.18$ and $17.7$, and is therefore
nonminimum phase. Moreover, it has a pole at the origin,
so the open-loop system is not asymptotically stable.

The pair $(A,B)$ is controllable and the pair $(A,C)$
is observable. For simulation, the initial condition is selected as
$x(0)=[1,1,1,1,1,1]^{\mathrm T}/\sqrt{6}$.
Throughout this example, $h=T/3$.

\subsubsection{Without Measurement Delays}

To investigate the effect of the prescribed time on
transient performance, we consider $T=60\,\mathrm{s}$
and $T=30\,\mathrm{s}$.
The weighting matrix is fixed as $Q_{\mathrm d}=I_6$.
For each prescribed time $T$, we apply the multiple-delay controller
\eqref{eq:multiple_delay_controller} in
Theorem~\ref{thm:multiple_delay_PT} with six artificial delays,
collected in the vector $\boldsymbol h_T=[h_1,\ldots,h_6]^{\mathrm T}$.
The superscript $\mathrm u$ denotes the uniformly spaced baseline,
with $h_i^{\mathrm u}=h+ih/7$ for $i=1,\ldots,6$ and $h=T/3$.
For $T=60\,\mathrm{s}$, this gives
\begin{equation}
\boldsymbol h_{60}^{\mathrm u}
=[22.8571,25.7143,28.5714,31.4286,34.2857,37.1429]^{\mathrm T},
\label{eq:uniform_delays_T60}
\end{equation}
whereas for $T=30\,\mathrm{s}$, we obtain
\[
\boldsymbol h_{30}^{\mathrm u}
=[11.4286,12.8571,14.2857,15.7143,17.1429,18.5714]^{\mathrm T}.
\]

To reduce state overshoot and peak control magnitude, we tune the
artificial delays using the Nelder--Mead simplex method. Let $\boldsymbol h_T^{\mathrm o}$ denote the resulting
delay vector for the prescribed time $T$. We obtain
\[
\begin{aligned}
\boldsymbol h_{60}^{\mathrm o}
&=[20.0000,20.3160,20.8246,39.3747,39.7730,39.9995]^{\mathrm T},\\
\boldsymbol h_{30}^{\mathrm o}
&=[10.0016,10.1017,10.7284,15.7086,19.8994,20.0000]^{\mathrm T}.
\end{aligned}
\]

The resulting trajectories are shown in
Fig.~\ref{figure2}.
For both prescribed times, the state converges to the
origin by $t=T$ and remains there thereafter.
Compared with the uniformly spaced delays, the optimized
delays reduce the transient peaks of both the state norm
and the control input.
Moreover, the shorter prescribed time requires larger
control peaks, whereas the longer waiting interval for
$T=60\,\mathrm{s}$ allows a larger state excursion before
convergence.

\subsubsection{Sampled-Data Delayed Implementation}

We next apply the sampled-data controller from
Corollary~\ref{rem:sampling_delay_rank_preservation}
to the same flexible-structure system, using the output samples
in \eqref{eq:sampled_output_measurements}.

We select $T=60\,\mathrm{s}$, $h=20\,\mathrm{s}$,
$d_{\mathrm s}=2/7\,\mathrm{s}$, and
$\tau_M=d_{\mathrm s}$.
The six values are chosen as
$\boldsymbol h_{60}^{\mathrm u}$ given by (\ref{eq:uniform_delays_T60}), using the exact values
$h_i=20+20i/7=(70+10i)d_{\mathrm s}$,
$i=1,\ldots,6$ and $\nu_i=70+10i$.
Thus, they are integer multiples of the sampling period
and satisfy $h\leq h_1<\cdots<h_6\leq2h-\tau_M$.
According to
Corollary~\ref{rem:sampling_delay_rank_preservation},
the condition of
Lemma~\ref{thm:known_tv_Gramian_positive}
is satisfied, and $G_{\mathrm v}>0$.
The weighting matrix is selected as $Q_{\mathrm v}=I_6$,
and the gain functions are constructed from
(\ref{eq:known_tv_gain}).

\begin{figure}[t]
\centering
\includegraphics[scale=0.55]{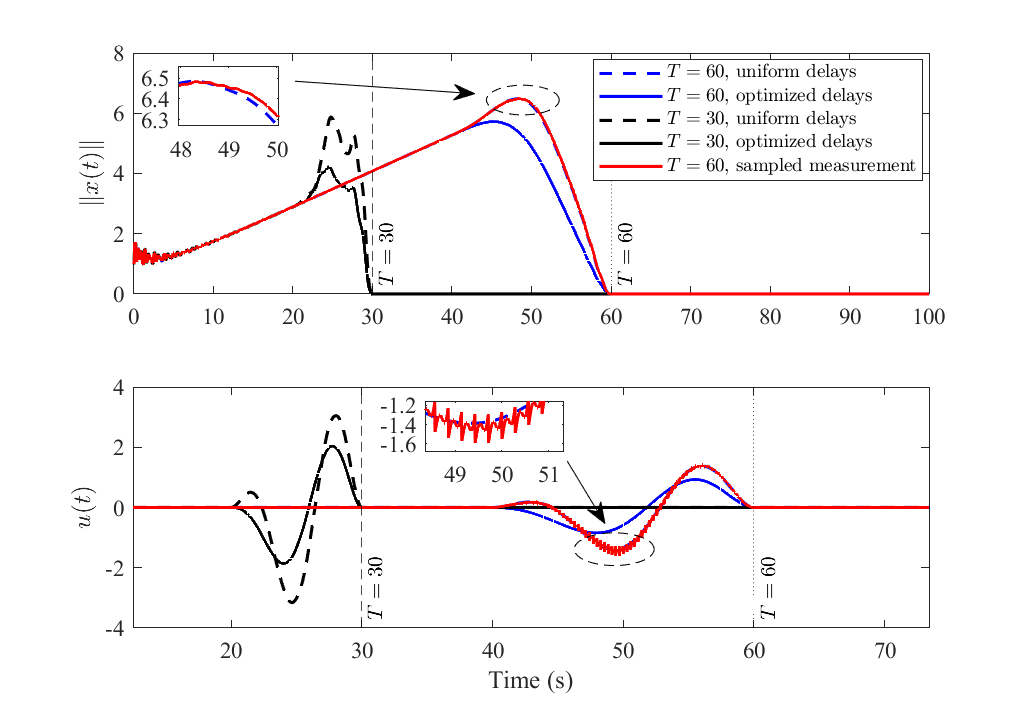}\caption{State and controller
trajectories of the flexible-structure system. }%
\label{figure2}%
\end{figure}

The red curves in Fig.~\ref{figure2} show the trajectories
under sampled-output measurements.
The state converges to the origin by $T=60\,\mathrm{s}$ and remains there
thereafter.
Compared with the continuous-output case using the same
uniformly spaced artificial delays, the sampled-output
case exhibits small oscillations and slightly larger
control peaks.
\subsection{Cart--Pole System With Position-Only Output}
\label{ex:cart}

To further demonstrate the applicability of the proposed
method to a physical nonminimum-phase system, consider
the linearized cart--pole system adopted from
\cite{CTMSpole}.
We consider only the cart-position output.
The system is described by (\ref{eq:MIMO_plant}) with
$x(t)=[x_{\mathrm c}(t),\dot{x}_{\mathrm c}(t),
\theta(t),\dot{\theta}(t)]^{\mathrm T}$ and
\begin{align*}
A&=
\begin{bmatrix}
0 & 1 & 0 & 0\\
0 & -\dfrac{(I+ml^2)b_{\mathrm c}}{p_0}
  & \dfrac{m^2gl^2}{p_0} & 0\\
0 & 0 & 0 & 1\\
0 & -\dfrac{mlb_{\mathrm c}}{p_0}
  & \dfrac{mgl(M+m)}{p_0} & 0
\end{bmatrix},\\
B&=
\begin{bmatrix}
0 & \dfrac{I+ml^2}{p_0} & 0 & \dfrac{ml}{p_0}
\end{bmatrix}^{\mathrm T},
\qquad
C=
\begin{bmatrix}
1&0&0&0
\end{bmatrix},
\end{align*}
where $p_0=I(M+m)+Mml^2$.
Here, $x_{\mathrm c}(t)$ denotes the cart displacement,
and $\theta(t)$ denotes the angular displacement of
the pendulum from the upright position.
Similar to \cite{CTMSpole}, the physical parameters
are selected as
$M=0.5\,\mathrm{kg}$, $m=0.2\,\mathrm{kg}$,
$b_{\mathrm c}=0.1\,\mathrm{N\,s/m}$,
$I=0.006\,\mathrm{kg\,m^2}$, $l=0.3\,\mathrm m$,
and $g=9.8\,\mathrm{m/s^2}$. We impose the control input constraint $|u(t)|\leq15\,\mathrm{N}$.

The transfer function from $u$ to $x_{\mathrm c}$ is
\[
G_{\mathrm p}(s)=
\frac{1.8182s^2-44.5455}
{s^4+0.1818s^3-31.1818s^2-4.4545s}.
\]
Its zeros are
\begin{equation}
s=\pm\sqrt{\frac{mgl}{I+ml^2}}=\pm4.9497.
\label{eq:cart_position_zeros}
\end{equation}
Thus, the system has a right-half-plane zero at
$s=4.9497$ and is nonminimum phase.

The pair $(A,B)$ is controllable, and
$(A,C)$ is observable.
We select the prescribed time $T=0.9\,\mathrm s$, so that
$h=T/3=0.3\,\mathrm s$.
Four artificial delays are used:
$h_1=0.3\,\mathrm s$, $h_2=0.4\,\mathrm s$,
$h_3=0.5\,\mathrm s$, and $h_4=0.6\,\mathrm s$.

For the controller in Theorem~\ref{thm:multiple_delay_PT}, we set $Q_{\mathrm d}=I_4$. The resulting Gramian $G$, defined in
\eqref{eq:multiple_delay_Gramian}, has a condition number of
approximately $1.17\times10^{14}$.
To avoid explicit Gramian inversion, we adopt the
numerical implementation in
Section~\ref{subsec:finite_dimensional_implementation},
using $N_b=8$ shifted Legendre basis functions.
The matrix $K_b$ defined in
(\ref{eq:finite_K_matrix}) has full row rank,
$\operatorname{rank}(K_b)=n^2=16$,
and a condition number of approximately
$5.69\times10^7$, substantially lower than that of $G$.
The solutions of the terminal cancellation equations in
\eqref{eq:finite_terminal_cancellation} can be written as
$a=K_b^\dagger b+N\eta$, where $b$ is the right-hand side of
$K_ba=b$ and the columns of $N$ span $\ker K_b$.
We use sequential quadratic programming to select $\eta$,
minimizing the worst-case peak control magnitude over the
initial-state ellipse.

\begin{figure}[t]
\centering
\includegraphics[scale=0.5]{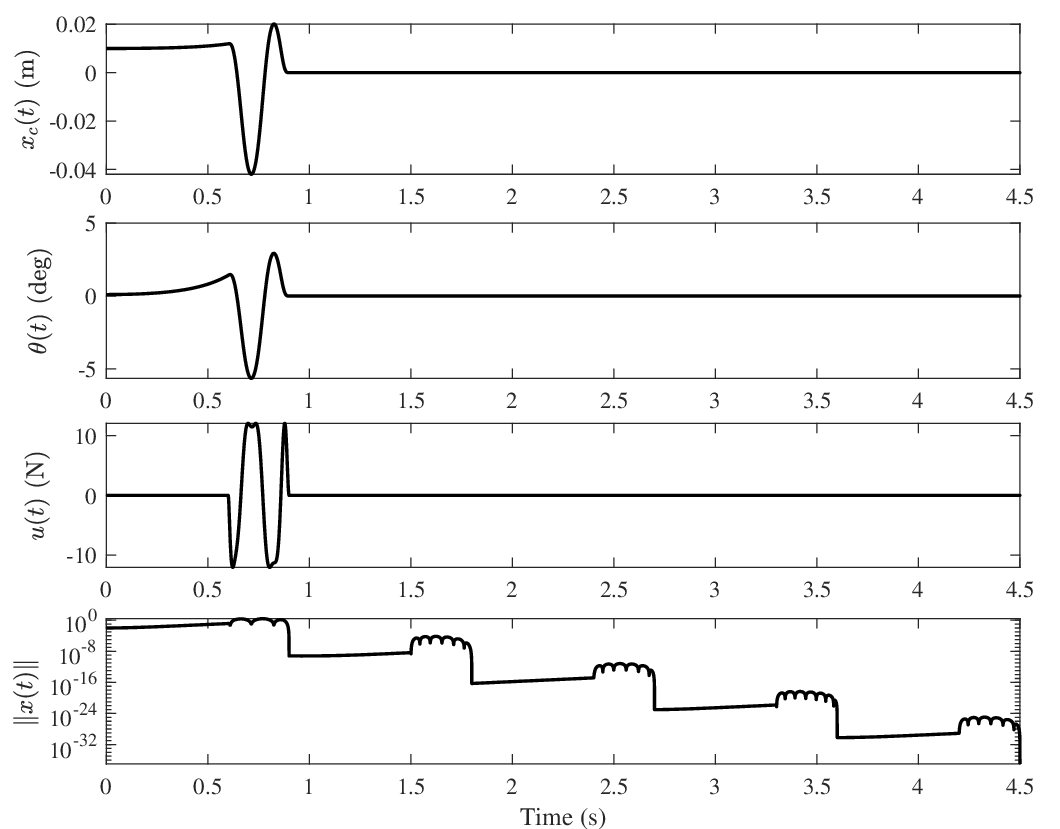}
\caption{State and control trajectories of the cart--pole
system with position-only output feedback.}
\label{fig:cart_pole_comparison}
\end{figure}

For simulation, take the initial condition
$x(0)=[0.01,0,0.1\pi/180,0]^{\mathrm T}$. The simulation results are shown in
Fig.~\ref{fig:cart_pole_comparison}.
The proposed controller drives the state to the origin
by the prescribed time $T=0.9\,\mathrm s$, using only cart-position
measurements. The peak control force in the simulation is approximately
$12.8\,\mathrm{N}$, below the specified limit of $15\,\mathrm{N}$.

\section{Conclusion}
\label{Conclu}
This paper has provided a constructive solution to the
static output-feedback stabilization problem and further
achieved prescribed-time stabilization for every controllable
and observable LTI system using two methods.
The first employs multiple artificial output delays
with smooth time-varying gains and stabilizes every
controllable and observable LTI system using a
generalized controllability--observability Gramian.
The second, simpler method combines constant-gain
feedback of the current output with time-varying-gain
feedback of a single artificially delayed output
and also applies to every controllable and observable LTI system.
The multiple-delay design has been extended to known
bounded and fast-varying (without any constraints on the delay derivative) measurement delays, yielding fixed-time
stabilization under a corresponding rank condition, as well as to
sampled-data control.
For unknown constant measurement or input delays,
fixed-time stabilization has been established for systems with
nilpotent state matrices.
A numerical implementation of Gramian-based gains
has also been suggested. Numerical examples involving nonminimum-phase,
flexible-structure, and cart--pole systems have illustrated
the effectiveness of the proposed methods.
Future work will address robustness to model uncertainty,
measurement noise, and delay uncertainty.

\section*{Appendix}
\subsection*{A1: The Proof of Lemma \ref{cor:equally_spaced_delays}}
Define $\mathit{\Phi}_{d}:=\mathrm{e}^{-Ad}$.
Since $h_{i}=h_{0}+(i-1)d$, we
have
$C\mathrm{e}^{-Ah_{i}}=
C\mathrm{e}^{-A[h_{0}+(i-1)d]}=
C\mathit{\Phi}_{d}^{\,i-1}\mathrm{e}^{-Ah_{0}}$.
Consequently, for $n_{0}=n$, (\ref{eq:sampled_observability_matrix}) can be rewritten as
\begin{align}
\mathcal{O}_{\boldsymbol h}
=
\begin{bmatrix}
C\mathrm{e}^{-Ah_{1}}\\
C\mathrm{e}^{-Ah_{2}}\\
\vdots\\
C\mathrm{e}^{-Ah_{n}}
\end{bmatrix}=
\begin{bmatrix}
C\\
C\mathit{\Phi}_{d}\\
\vdots\\
C\mathit{\Phi}_{d}^{\,n-1}
\end{bmatrix}
\mathrm{e}^{-Ah_{0}}.
\label{eq:sampled_observability_factorization}
\end{align}
Since $\mathrm{e}^{-Ah_{0}}$ is nonsingular,
\[
\operatorname{rank}\mathcal O_{\boldsymbol h}
=
\operatorname{rank}\bigl(
\operatorname{col}(C,C\mathit{\Phi}_d,\ldots,
C\mathit{\Phi}_d^{\,n-1})
\bigr).
\]
Therefore, it is sufficient to prove that the discrete-time pair
$(\mathit{\Phi}_{d},C)$ is observable for every sufficiently small $d>0$.

Let $\lambda_1,\ldots,\lambda_{\varsigma}$ denote
the distinct eigenvalues of $A$, where
$\varsigma\leq n$ is the number of distinct eigenvalues. The corresponding eigenvalues of $\mathit{\Phi}_{d}$ are
$\mathrm{e}^{-\lambda_{1}d},\ldots,
\mathrm{e}^{-\lambda_{\varsigma}d}$.
Two distinct eigenvalues $\lambda_{j}$ and $\lambda_{k}$ of $A$ are mapped to
the same eigenvalue of $\mathit{\Phi}_{d}$ only if $\mathrm{e}^{-\lambda_{j}d}
=\mathrm{e}^{-\lambda_{k}d}$
or equivalently, $(\lambda_{j}-\lambda_{k})d
=
2\pi\mathrm{i}\ell$
for some nonzero integer $\ell$.

Taking absolute values gives
$
|\lambda_j-\lambda_k|d=2\pi|\ell|\geq2\pi$.
Thus, when $\varsigma\geq2$, define
\[
d_A:=
\frac{2\pi}{\max_{j\neq k}|\lambda_j-\lambda_k|}.
\]
For $\varsigma=1$, set $d_A=+\infty$.
Then
\begin{equation}
\mathrm e^{-\lambda_jd}\neq\mathrm e^{-\lambda_kd},
\qquad j\neq k,\quad 0<d<d_A.
\label{eq:no_sampling_aliasing}
\end{equation}

Under condition (\ref{eq:no_sampling_aliasing}),
distinct eigenvalues of $A$ remain distinct under
the exponential mapping.
Moreover, since $-d\mathrm e^{-\lambda d}\neq0$,
each Jordan block $J$ and $\mathrm e^{-dJ}$
have the same eigenvectors.
Therefore, every nonzero eigenvector $v$ of
$\mathit{\Phi}_d$ is also an eigenvector of $A$. In particular,
there exists an eigenvalue $\lambda$ of $A$ such that
$Av=\lambda v$.
Since $(A,C)$ is observable, the PBH observability criterion gives
$Cv\neq0$
for every nonzero eigenvector $v$ of $A$. It follows that no nonzero
eigenvector of $\mathit{\Phi}_{d}$ belongs to the null space of $C$. By the PBH
criterion again, $(\mathit{\Phi}_{d},C)$ is observable for every
$0<d<d_{A}$.

Finally, since $h\leq h_0<2h$, choose
\begin{equation}
0<d^\ast\leq
\min\left\{
d_A,\frac{2h-h_0}{n-1}
\right\},
\label{eq:admissible_delay_spacing}
\end{equation}
when $n\geq2$. Then, for every $0<d<d^{\ast}$, all artificial delays satisfy
$h\leq h_{1}<h_{2}<\cdots<h_{n}<2h$,
and the pair $(\mathit{\Phi}_{d},C)$ is observable. Hence, its $n$-step observability
matrix has full column rank:
$\operatorname{rank}\bigl(\operatorname{col}(C,C\mathit{\Phi}_d,\ldots,C\mathit{\Phi}_d^{\,n-1})\bigr)=n$.
Using (\ref{eq:sampled_observability_factorization}), we conclude that
$\operatorname{rank}\mathcal{O}_{\boldsymbol h}=n$.
For $n=1$, the conclusion follows directly from the observability of $(A,C)$.
\subsection*{A2: The Proof of Lemma \ref{thm:multiple_delay_Gramian}}
For $s\in[2h,3h]$, define
$\mathcal{L}(s) := \mathrm{e}^{A(3h-s)}B \in \mathbb{R}^{n\times m}$
and
$
\mathcal{R}_{i}(s) := C\mathrm{e}^{A(s-h_{i}-2h)} \in \mathbb{R}^{p\times n}$,
$i=1,\ldots,n_{0}$.
Then the $i$-th block column of
$\mathit{\Gamma}_{\mathrm d}(s)$ can be written as
\begin{align*}
\mathit{\Gamma}_{\mathrm{d},i}(s)  & = \mathrm{e}^{A^{\mathrm{T}}(s-h_{i}%
-2h)}C^{\mathrm{T}} \otimes \mathrm{e}^{A(3h-s)}B\\
& = \mathcal{R}_{i}^{\mathrm{T}}(s) \otimes \mathcal{L}(s) \in \mathbb{R}%
^{n^{2}\times mp}.
\end{align*}
Consequently,
\[
\mathit{\Gamma}_{\mathrm{d}}(s) =
\begin{bmatrix}
\mathcal{R}_{1}^{\mathrm{T}}(s)\otimes \mathcal{L}(s) & \cdots & \mathcal{R}%
_{n_{0}}^{\mathrm{T}}(s) \otimes \mathcal{L}(s)
\end{bmatrix}
.
\]

Collect all the delayed-output matrices in
\begin{align*}
\mathcal{R}(s):=  &
\begin{bmatrix}
\mathcal{R}_{1}(s)\\
\vdots \\
\mathcal{R}_{n_{0}}(s)
\end{bmatrix}
=%
\begin{bmatrix}
C\mathrm{e}^{A(s-h_{1}-2h)}\\
\vdots \\
C\mathrm{e}^{A(s-h_{n_{0}}-2h)}%
\end{bmatrix}
\\
=  & \mathcal{O}_{\boldsymbol{h}}\mathrm{e}^{A(s-2h)}\in \mathbb{R}%
^{n_{0}p\times n}.
\end{align*}
By condition (\ref{eq:sampled_observability_rank}), $\operatorname{rank}%
\! \left(  \mathcal{O}_{\boldsymbol{h}}\right)  =n$. Moreover, $\mathrm{e}%
^{A(s-2h)}$ is nonsingular for every $s$. Therefore,
\[
\operatorname{rank}\! \left(  \mathcal{R}(s)\right)  =\operatorname{rank}%
\! \left(  \mathcal{O}_{\boldsymbol{h}}\mathrm{e}^{A(s-2h)}\right)  =n,\qquad
s\in \lbrack2h,3h].
\]
Thus, $\mathcal{R}(s)$ has full column rank, or equivalently, $\mathcal{R}%
^{\mathrm{T}}(s)$ has full row rank. Hence, for every $s\in \lbrack2h,3h]$,
there exists a right inverse $\mathcal J(s)\in \mathbb{R}%
^{n_{0}p\times n}$ such that
$\mathcal{R}^{\mathrm{T}}(s)\mathcal J(s)=I_{n}$.

We next characterize the null space of $\mathit{\Gamma}_{\mathrm{d}%
}^{\mathrm{T}}(s)$. Let $z\in \mathbb{R}^{n^{2}}$ and define
$
Z:=\operatorname{vec}^{-1}(z)\in \mathbb{R}^{n\times n}$.
Suppose that
$\mathit{\Gamma}_{\mathrm{d}}^{\mathrm{T}}(s)z=0$, $s\in(2h,3h)$.
Since $\mathit{\Gamma}_{\mathrm{d}}(s)$ is composed of the blocks
$\mathit{\Gamma}_{\mathrm{d},i}(s)$, the preceding equality implies
\[
\mathit{\Gamma}_{\mathrm{d},i}^{\mathrm{T}}(s)z=0,\qquad i=1,\ldots
,n_{0}.
\]
Using
$\mathit{\Gamma}_{\mathrm{d},i}(s)=\mathcal{R}_{i}^{\mathrm{T}}(s)\otimes
\mathcal{L}(s)$,
we obtain
\[
\mathit{\Gamma}_{\mathrm{d},i}^{\mathrm{T}}(s)z=\left(  \mathcal{R}%
_{i}(s)\otimes \mathcal{L}^{\mathrm{T}}(s)\right)  \operatorname{vec}(Z).
\]
Applying the Kronecker--vectorization identity
$(B^{\mathrm{T}}\otimes A)\operatorname{vec}(X)=\operatorname{vec}(AXB)$
with $A=\mathcal{L}^{\mathrm{T}}(s)$, $X=Z$, and $B=\mathcal{R}_{i}^{\mathrm{T}%
}(s)$ gives
\[
\mathit{\Gamma}_{\mathrm{d},i}^{\mathrm{T}}(s)z=\operatorname{vec}\! \left(
\mathcal{L}^{\mathrm{T}}(s)Z\mathcal{R}_{i}^{\mathrm{T}}(s)\right)  .
\]
Therefore, $\mathit{\Gamma}_{\mathrm{d},i}^{\mathrm{T}}(s)z=0$ is equivalent
to
\[
\mathcal{L}^{\mathrm{T}}(s)Z\mathcal{R}_{i}^{\mathrm{T}}(s)=0,\qquad
i=1,\ldots,n_{0}.
\]
Collecting these equalities for all $i$ yields $\mathcal{L}^{\mathrm{T}%
}(s)Z\mathcal{R}^{\mathrm{T}}(s)=0$. Multiplying this equality from the right
by $\mathcal J(s)$ gives
\begin{align*}
\mathcal{L}^{\mathrm{T}}(s)Z  =\mathcal{L}^{\mathrm{T}}(s)Z\mathcal{R}%
^{\mathrm{T}}(s)\mathcal J(s)  =0,\qquad s\in(2h,3h).
\end{align*}

It remains to use the controllability of $(A,B)$. We claim that
\begin{equation}
\operatorname{span} \left \{  \operatorname{Im}\mathcal{L}(s): s\in(2h,3h)
\right \}  = \mathbb{R}^{n}.\label{eq:L_span}%
\end{equation}
Suppose, to the contrary, that a nonzero vector
$v\in\mathbb R^n$ is orthogonal to this span. Then
\[
v^{\mathrm{T}}\mathcal{L}(s) = v^{\mathrm{T}}\mathrm{e}^{A(3h-s)}B = 0, \qquad
s\in(2h,3h).
\]
The left-hand side is an analytic matrix-valued function of $s$. Fix an
arbitrary $s_{0}\in(2h,3h)$. Since this function vanishes on an open interval,
all of its derivatives at $s=s_{0}$ vanish. In particular,
\[
\frac{\mathrm{d}^{k}}{\mathrm{d}s^{k}} \mathrm{e}^{A(3h-s)} = (-1)^{k}
\mathrm{e}^{A(3h-s)}A^{k},
\]
and hence
\[
(-1)^{k} v^{\mathrm{T}}\mathrm{e}^{A(3h-s_{0})}A^{k}B = 0, \qquad
k=0,1,\ldots,n-1.
\]
Define
\[
\widetilde{v}^{\mathrm{T}} := v^{\mathrm{T}}\mathrm{e}^{A(3h-s_{0})}.
\]
Then
\[
\widetilde{v}^{\mathrm{T}}
\begin{bmatrix}
B & AB & \cdots & A^{n-1}B
\end{bmatrix}
= 0.
\]
Since $(A,B)$ is controllable,
\[
\operatorname{rank}
\begin{bmatrix}
B & AB & \cdots & A^{n-1}B
\end{bmatrix}
= n.
\]
Therefore, $\widetilde{v}=0$. Since $\mathrm{e}^{A(3h-s_{0})}$ is nonsingular,
this implies $v=0$. Thus, the orthogonal complement of the span in
(\ref{eq:L_span}) is trivial, which proves (\ref{eq:L_span}).

Recall that
$
\mathcal{L}^{\mathrm{T}}(s)Z=0$ for all $ s\in(2h,3h)$.
This means that $x^{\mathrm{T}}Z=0$ for every
\[
x\in \operatorname{span}\left \{  \operatorname{Im}\mathcal{L}(s):s\in
(2h,3h)\right \}  .
\]
By (\ref{eq:L_span}), this span is equal to $\mathbb{R}^{n}$. Consequently,
\[
x^{\mathrm{T}}Z=0,\qquad \forall \,x\in \mathbb{R}^{n},
\]
which implies $Z=0$. Since $z=\operatorname{vec}(Z)$, we also have $z=0$. We
have thus proved that
\begin{equation}
\left[  \mathit{\Gamma}_{\mathrm{d}}^{\mathrm{T}}(s)z=0\  \text{for every }%
s\in(2h,3h)\right]  \quad \Longrightarrow \quad
z=0.\label{eq:Gamma_trivial_common_kernel}%
\end{equation}

We now prove the positive definiteness of $G$. For any
$z\in \mathbb{R}^{n^{2}}$, it follows from (\ref{eq:multiple_delay_Gramian})
that
\begin{align*}
z^{\mathrm{T}}Gz ={} &  \int_{2h}^{3h} S(s) z^{\mathrm{T}%
}\mathit{\Gamma}_{\mathrm{d}}(s) Q_{\mathrm{d}} \mathit{\Gamma}_{\mathrm{d}%
}^{\mathrm{T}}(s)z \mathrm{d}s\\
\geq{} &  \lambda_{\min}(Q_{\mathrm{d}}) \int_{2h}^{3h} S(s) \left\|
\mathit{\Gamma}_{\mathrm{d}}^{\mathrm{T}}(s)z \right\| ^{2}  \mathrm{d}s
\geq  0.
\end{align*}
Thus, $G$ is positive semidefinite.

Suppose that $z^{\mathrm{T}}Gz=0$. Since $Q_{\mathrm{d}}>0$,
$S(s)>0$ on $(2h,3h)$, and
$\mathit{\Gamma}_{\mathrm{d}}^{\mathrm{T}}(s)z$ is continuous in $s$, the
preceding integral equality implies
\[
\mathit{\Gamma}_{\mathrm{d}}^{\mathrm{T}}(s)z=0,
\qquad s\in(2h,3h).
\]
It then follows from (\ref{eq:Gamma_trivial_common_kernel}) that $z=0$.
Therefore,
\[
z^{\mathrm{T}}Gz>0,\qquad \forall \,z\in \mathbb{R}^{n^{2}}%
\setminus \{0\},
\]
and hence $G>0$.

\subsection*{A3: The Proof of Lemma \ref{thm:MIMO_constant_Gramian}}

The eigenvalue arguments are carried out over $\mathbb C$.
Since the matrices in the controllability condition are real,
their ranks are the same over $\mathbb R$ and $\mathbb C$.
Let $\lambda_1,\ldots,\lambda_{\varsigma}$ denote the
distinct eigenvalues of $A$.

\medskip
\noindent\textbf{Step 1: Selection of $H$, $v$, and $w$.}
By the output-feedback result in \cite{davison73tac},
Assumption~\ref{assumption_minimal} guarantees the existence
of a real matrix $H\in\mathbb R^{m\times p}$ such that
$F_H=A+BHC$ has $n$ distinct eigenvalues
$\zeta_1,\ldots,\zeta_n$.
State feedback preserves controllability, and output
injection preserves observability. Therefore,
$(F_H,B)$ is controllable and $(F_H,C)$ is observable.

Let $r_i$ and $l_i^{\mathrm T}$ be nonzero right and left
eigenvectors of $F_H$ associated with $\zeta_i$.
The PBH criterion gives
\[
l_i^{\mathrm T}B\neq0,\qquad Cr_i\neq0,
\qquad i=1,\ldots,n.
\]
Consequently, the real vectors $v\in\mathbb R^m$ for which
$l_i^{\mathrm T}Bv=0$ form a proper real linear subspace
for each $i$. A finite union of proper linear subspaces
cannot cover $\mathbb R^m$. Hence, a real vector $v$ can
be selected so that $l_i^{\mathrm T}Bv\neq0$ for every $i$.
Similarly, a real vector $w\in\mathbb R^p$ can be selected
so that $w^{\mathrm T}Cr_i\neq0$ for every $i$.
Writing
\[
b=Bv\in\mathbb R^n,\qquad
c=w^{\mathrm T}C\in\mathbb R^{1\times n},
\]
the PBH criterion shows that $(F_H,b)$ is controllable
and $(F_H,c)$ is observable.

For $K_\alpha=H+\alpha vw^{\mathrm T}$, we have
\begin{equation}
F_{K_\alpha}=F_H+\alpha bc.
\label{eq:rank_one_closed_loop_family}
\end{equation}
Thus, a rank-one output-feedback family can be used
after the preliminary selection of $H$.

\medskip
\noindent\textbf{Step 2: Eigenvalue conditions outside a finite set.}
Define
\[
D(z)=\det(zI_n-F_H),\qquad
U(z)=c\operatorname{adj}(zI_n-F_H)b.
\]
Since $(F_H,b,c)$ is a minimal SISO realization, $D$ and
$U$ are coprime, $U\not\equiv0$, $\deg D=n$, and
$\deg U\leq n-1$. Using (\ref{eq:rank_one_closed_loop_family}) and the
rank-one determinant identity
$\det(M+uv^{\mathrm T})
=\det(M)+v^{\mathrm T}\operatorname{adj}(M)u$,
where $M\in\mathbb C^{n\times n}$ and $u,v\in\mathbb C^n$,
with $M=zI_n-F_H$, $u=-\alpha b$, and $v=c^{\mathrm T}$,
we obtain
\[
P_\alpha(z):=\det(zI_n-F_{K_\alpha})
=D(z)-\alpha U(z).
\]

\smallskip
\noindent\textit{(i) Distinct eigenvalues of $F_{K_\alpha}$}.
If $z$ is a repeated root of $P_\alpha$, then
\[
D(z)-\alpha U(z)=0,\qquad
D'(z)-\alpha U'(z)=0.
\]
Coprimeness implies $U(z)\neq0$, and hence
\[
\alpha=\frac{D(z)}{U(z)},\qquad
V(z):=D'(z)U(z)-D(z)U'(z)=0.
\]
The polynomial $V$ is not identically zero. Otherwise,
$(D/U)'\equiv0$, so the rational function $D/U$ would
be constant, contradicting $\deg D>\deg U$.
Therefore, the set
\[
\mathcal Z_0:=
\left\{\frac{D(z)}{U(z)}:
z\in\mathbb C,\ V(z)=0,\ U(z)\neq0\right\}
\cap\mathbb R
\]
is finite, and $F_{K_\alpha}$ has $n$ distinct eigenvalues
whenever $\alpha\notin\mathcal Z_0$.

\smallskip
\noindent\textit{(ii) Exclusion of eigenvalue-difference coincidences.}
For $\alpha\notin\mathcal Z_0$, let
$\mu_1,\ldots,\mu_n$ denote the distinct eigenvalues
of $F_{K_\alpha}$.
Since
$A_{K_\alpha}=A^{\mathrm T}\otimes I_n
-I_n\otimes F_{K_\alpha}$,
its eigenvalues are the differences $\lambda_j-\mu_i$.
We next exclude coincidences between these differences
for distinct index pairs $(j,i)$.
For $j\neq k$, set $d_{jk}=\lambda_k-\lambda_j\neq0$.
If
$\lambda_j-\mu_i=\lambda_k-\mu_\ell$,
then $\mu_\ell=\mu_i+d_{jk}$.
Since both are roots of $P_\alpha=D-\alpha U$, we can get
\[
D(\mu_i)=\alpha U(\mu_i),\qquad
D(\mu_i+d_{jk})=\alpha U(\mu_i+d_{jk}).
\]
Eliminating $\alpha$ gives
\[
D(\mu_i)U(\mu_i+d_{jk})
-D(\mu_i+d_{jk})U(\mu_i)=0.
\]
This motivates defining
\[
R_{jk}(z)
:=D(z)U(z+d_{jk})-D(z+d_{jk})U(z).
\]
We next show that $R_{jk}$ is not identically zero.
Suppose, to the contrary, that $R_{jk}\equiv0$. Then
\[
D(z)U(z+d_{jk})=D(z+d_{jk})U(z).
\]
If $U(z)U(z+d_{jk})=0$, then $z$ and $z+d_{jk}$
cannot both be roots of $P_\alpha$, since $D$ and $U$
are coprime. At every point where $U(z)U(z+d_{jk})\neq0$,
division gives
\[
\frac{D(z+d_{jk})}{U(z+d_{jk})}
=\frac{D(z)}{U(z)}.
\]
Since $U$ is a nonzero polynomial, the excluded points
form a finite set, and the preceding equality holds
as an identity of rational functions.
Thus, $f=D/U$ has the nonzero period $d_{jk}$.

If $f$ has a pole at $z_0$, periodicity implies that
$z_0+\ell d_{jk}$ is also a pole for every
$\ell\in\mathbb Z_0^+$.
Since $d_{jk}\neq0$, these points are all distinct,
giving infinitely many poles.
This is impossible because a rational function has
only finitely many finite poles.
If it has no poles, it is a polynomial, and a nonconstant
polynomial cannot have a nonzero period. Thus, in either case, the nonconstant rational function
$D/U$ cannot have a nonzero period, and hence
$R_{jk}\not\equiv0$.

As shown above, any coincidence
$\lambda_j-\mu_i=\lambda_k-\mu_\ell$ with $j\neq k$
implies $R_{jk}(\mu_i)=0$.
Moreover, $U(\mu_i)\neq0$: otherwise,
$P_\alpha(\mu_i)=0$ would imply $D(\mu_i)=0$,
contradicting the coprimeness of $D$ and $U$.
The same argument gives
$U(\mu_i+d_{jk})=U(\mu_\ell)\neq0$.
Hence,
$\alpha=D(\mu_i)/U(\mu_i)$.

Since $R_{jk}$ is a nonzero polynomial, it has only
finitely many roots. Consequently, every real $\alpha$
producing such a coincidence belongs to the finite set
\[
\mathcal Z_{jk}:=
\left\{\frac{D(z)}{U(z)}:
z\in\mathbb C, R_{jk}(z)=0,
U(z)U(z+d_{jk})\neq0\right\}\cap\mathbb R.
\]
Define
\[
\mathcal Z=\mathcal Z_0\cup
\bigcup_{\substack{1\leq j,k\leq\varsigma\\j\neq k}}
\mathcal Z_{jk}.
\]
This set is finite because it is a finite union of
finite sets.
For every $\alpha\in\mathbb R\setminus\mathcal Z$,
$F_{K_\alpha}$ has $n$ distinct eigenvalues
$\mu_1,\ldots,\mu_n$, and no coincidence is possible
when $j\neq k$.
If $j=k$, equality of the differences implies
$\mu_i=\mu_\ell$, and hence $i=\ell$.
Therefore,
\begin{equation}
\lambda_j-\mu_i=\lambda_k-\mu_\ell
\quad\Longrightarrow\quad(j,i)=(k,\ell).
\label{eq:MIMO_difference_separation}
\end{equation}

\medskip
\noindent\textbf{Step 3: Controllability through a Sylvester equation.}
Fix $\alpha\notin\mathcal Z$ and write $K=K_\alpha$.
By the PBH controllability criterion, it is sufficient
to show that, for every $\eta\in\mathbb C$,
\[
A_K^{\mathrm T}z=\eta z,\qquad B_C^{\mathrm T}z=0
\]
imply $z=0$. Write $z=\operatorname{vec}(Z)$, where
$Z\in\mathbb C^{n\times n}$. Since
\[
A_K^{\mathrm T}
=A\otimes I_n-I_n\otimes F_K^{\mathrm T},
\qquad
B_C^{\mathrm T}=C\otimes B^{\mathrm T},
\]
writing $z=\operatorname{vec}(Z)$ and using the identity
$\operatorname{vec}(MXN)
=(N^{\mathrm T}\otimes M)\operatorname{vec}(X)$
gives
\begin{align*}
A_K^{\mathrm T}z
&=(A\otimes I_n)\operatorname{vec}(Z)
 -(I_n\otimes F_K^{\mathrm T})\operatorname{vec}(Z)\\
&=\operatorname{vec}(ZA^{\mathrm T})
 -\operatorname{vec}(F_K^{\mathrm T}Z)=\operatorname{vec}(ZA^{\mathrm T}-F_K^{\mathrm T}Z),
\end{align*}
and
\[
B_C^{\mathrm T}z
=(C\otimes B^{\mathrm T})\operatorname{vec}(Z)
=\operatorname{vec}(B^{\mathrm T}ZC^{\mathrm T}).
\]
Since the vectorization operator is injective,
the conditions $A_K^{\mathrm T}z=\eta z$ and
$B_C^{\mathrm T}z=0$ are equivalent to
\begin{equation}
ZA^{\mathrm T}-F_K^{\mathrm T}Z=\eta Z,
\qquad
B^{\mathrm T}ZC^{\mathrm T}=0.
\label{eq:MIMO_Sylvester_PBH}
\end{equation}

Since $F_K$ has $n$ distinct eigenvalues, there exists a
nonsingular matrix $R\in\mathbb C^{n\times n}$ such that
\begin{equation}
F_K^{\mathrm T}R
=R\mathit{\Lambda},
\label{eq:FK_transpose_diagonalization}
\end{equation}
where $\mathit{\Lambda}=\operatorname{diag}(\mu_1,\ldots,\mu_n)$. Set $Z=RY$, where $Y\in\mathbb C^{n\times n}$.
Substituting this expression into the Sylvester equation
in (\ref{eq:MIMO_Sylvester_PBH}) and using
(\ref{eq:FK_transpose_diagonalization}) gives
\[
YA^{\mathrm T}-\mathit{\Lambda} Y=\eta Y.
\]
Let $y_i\in\mathbb C^{1\times n}$ denote the $i$th row
of $Y$. Since the $i$th row of $\mathit{\Lambda} Y$ is $\mu_i y_i$,
the preceding equation implies
\begin{equation}
y_iA^{\mathrm T}=(\eta+\mu_i)y_i,
\qquad i=1,\ldots,n.
\label{eq:Sylvester_row_eigenvector}
\end{equation}
If $y_i\neq0$, then $y_i^{\mathrm T}$ is an eigenvector
of $A$ and $\eta=\lambda_j-\mu_i$ for some $j$.
Condition (\ref{eq:MIMO_difference_separation}) therefore
implies that at most one row of $Y$ is nonzero.

Suppose that $Z\neq0$.
The preceding argument shows that $Y$ has exactly one
nonzero row, say $y_i$.
Writing $R=[r_1,\ldots,r_n]$, we therefore obtain $Z=RY=r_i y_i$,
where $r_i$ is the $i$th column of $R$.
By (\ref{eq:FK_transpose_diagonalization}),
$r_i$ is an eigenvector of $F_K^{\mathrm T}$.
Since $(F_K,B)$ is controllable, the PBH criterion gives
$B^{\mathrm T}r_i\neq0$.
Similarly, (\ref{eq:Sylvester_row_eigenvector}) shows
that $y_i^{\mathrm T}$ is an eigenvector of $A$.
Observability of $(A,C)$ and the PBH criterion therefore
give $Cy_i^{\mathrm T}\neq0$, and hence
$y_iC^{\mathrm T}\neq0$.
It follows that
\[
B^{\mathrm T}ZC^{\mathrm T}
=(B^{\mathrm T}r_i)(y_iC^{\mathrm T})\neq0,
\]
contradicting (\ref{eq:MIMO_Sylvester_PBH}). Thus $Z=0$,
and the pair $(A_K,B_C)$ is controllable. In particular,
$\operatorname{rank}\mathcal{K}_{K_\alpha}=n^2$ for every
$\alpha\notin\mathcal Z$. Therefore,
$\mathcal E_{H,v,w}\subseteq\mathcal Z$ is finite.

\medskip
\noindent\textbf{Step 4: Positive definiteness of the Gramian.}
For any constant matrix $K$ satisfying
$\operatorname{rank}\mathcal{K}_{K}=n^2$ and any
$z\in\mathbb R^{n^2}$, we have
\[
z^{\mathrm T}G_Kz\geq\lambda_{\min}(Q)
\int_{2h}^{3h}S(s)
\|\mathit{\Gamma}_K^{\mathrm T}(s)z\|^2\,\mathrm ds\geq0.
\]
If $z^{\mathrm T}G_Kz=0$, continuity and positivity of
$S$ on $(2h,3h)$ imply
\[
z^{\mathrm T}\mathrm e^{A_K(s-3h)}B_C=0,
\qquad s\in(2h,3h).
\]
Fix $s_0\in(2h,3h)$. Taking successive derivatives gives
\[
z^{\mathrm T}\mathrm e^{A_K(s_0-3h)}A_K^jB_C=0,
\qquad j=0,\ldots,n^2-1.
\]
Thus $z^{\mathrm T}\mathrm e^{A_K(s_0-3h)}\mathcal{K}_{K}=0$.
Full row rank of $\mathcal{K}_{K}$ and nonsingularity
of the matrix exponential imply $z=0$. Hence $G_{K}>0$.

\subsection*{A4: Proof of Lemma \ref{thm:known_tv_Gramian_positive}}

Let $\mathcal{L}_{\mathrm v}(s)=\mathrm e^{A(3h-s)}B$ and
let $\mathcal{R}_{\mathrm v,i}(s)$ be the $i$th block row of
$
\mathcal{R}_{\mathrm v}(s)
:=\mathcal{O}_{\boldsymbol h,\tau}(s)\mathrm e^{A(s-2h)}$.
Then the $i$-th block column of
$\mathit{\Gamma}_{\mathrm v}(s)$ can be written as
$\mathit{\Gamma}_{\mathrm v,i}(s)
=\mathcal{R}_{\mathrm v,i}^{\mathrm T}(s)
\otimes\mathcal{L}_{\mathrm v}(s)$.
By (\ref{eq:effective_observability_rank}) and the
nonsingularity of the matrix exponential,
$\mathcal{R}_{\mathrm v}(s)$ has full column rank on
$\mathcal I$.

For any $z\in\mathbb R^{n^2}$, we have
\[
z^{\mathrm T}G_{\mathrm v}z
\geq\lambda_{\min}(Q_{\mathrm v})
\int_{\mathcal I}S(s)
\|\mathit{\Gamma}_{\mathrm v}^{\mathrm T}(s)z\|^2\,\mathrm ds
\geq0.
\]
Thus, $G_{\mathrm v}$ is positive semidefinite.
Suppose that $z^{\mathrm T}G_{\mathrm v}z=0$ and set
$Z=\operatorname{vec}^{-1}(z)$.
Continuity of $\mathit{\Gamma}_{\mathrm v}$ and positivity
of $S$ on $\mathcal I$ imply
$\mathit{\Gamma}_{\mathrm v}^{\mathrm T}(s)z=0$ throughout
$\mathcal I$.
As in the proof of Lemma~\ref{thm:multiple_delay_Gramian},
the Kronecker--vectorization identity gives
\[
0=\mathit{\Gamma}_{\mathrm v,i}^{\mathrm T}(s)z
=\operatorname{vec}\!\left(
\mathcal{L}_{\mathrm v}^{\mathrm T}(s)Z
\mathcal{R}_{\mathrm v,i}^{\mathrm T}(s)\right),
\quad i=1,\ldots,n_{0}.
\]
Collecting these equalities yields
$\mathcal{L}_{\mathrm v}^{\mathrm T}(s)Z
\mathcal{R}_{\mathrm v}^{\mathrm T}(s)=0$.
Since $\mathcal{R}_{\mathrm v}^{\mathrm T}(s)$ has a right
inverse, we obtain
$\mathcal{L}_{\mathrm v}^{\mathrm T}(s)Z=0$ on $\mathcal I$.

Next, let $v\in\mathbb R^n$ satisfy
$v^{\mathrm T}\mathcal{L}_{\mathrm v}(s)=0$ on $\mathcal I$.
Fix $s_0\in\mathcal I$. Differentiating at $s_0$ gives
\[
v^{\mathrm T}\mathrm e^{A(3h-s_0)}A^kB=0,
\qquad k=0,\ldots,n-1.
\]
Controllability of $(A,B)$ and nonsingularity of the
matrix exponential imply $v=0$. Hence the columns of
$\mathcal{L}_{\mathrm v}(s)$ over $s\in\mathcal I$ span
$\mathbb R^n$.
Since $\mathcal{L}_{\mathrm v}^{\mathrm T}(s)Z=0$,
each column of $Z$ is orthogonal to this span.
Thus $Z=0$ and $z=0$, proving $G_{\mathrm v}>0$.



\end{spacing}

\begin{thebibliography}{99}

\bibitem{Syrmos97auto}
V. L. Syrmos, C. T. Abdallah, P. Dorato, and K. Grigoriadis,
``Static output feedback---a survey,''
\textit{Automatica}, 33(2), 125--137, 1997.

\bibitem{sadabadi16}
M. S. Sadabadi and D. Peaucelle,
``From static output feedback to structured robust static output feedback: A survey,''
\textit{Annual Reviews in Control}, 42, 11--26, 2016.

\bibitem{lofberg26}
J. L\"ofberg,
``Static output-feedback stabilization is NP-hard,''
arXiv preprint arXiv:2609.16886, 2026.

\bibitem{barkai26}
G. Barkai and I. Shames,
``On the NP-hardness of unconstrained static output feedback stabilization,''
arXiv preprint arXiv:2609.38728, 2026.

\bibitem{ahmadi26}
A. A. Ahmadi, A. Chaudhry, I. Narang, and Y. Tang,
``Complexity of output feedback stabilization,''
arXiv preprint arXiv:2609.20636, 2026.

\bibitem{Brockett99}
R. W. Brockett, ``A stabilization problem,'' in
\textit{Open Problems in Mathematical Systems and Control Theory},
V. D. Blondel, E. D. Sontag, M. Vidyasagar, and J. C. Willems,
Eds. London, U.K.: Springer, pp. 75--78, 1999.

\bibitem{Moreau99cdc}
L. Moreau and D. Aeyels,
``Stabilization by means of periodic output feedback,'' in
\textit{Proceedings of the 38th IEEE Conference on Decision
and Control}, vol. 1, pp. 108--109, 1999.

\bibitem{Moreau04scl}
L. Moreau and D. Aeyels, ``Periodic output feedback
stabilization of single-input single-output continuous-time
systems with odd relative degree,''
\textit{Systems \& Control Letters}, 51(5), 395--406, 2004.

\bibitem{Allwright05auto}
J. C. Allwright, A. Astolfi, and H. P. Wong,
``A note on asymptotic stabilization of linear systems by
periodic, piecewise constant, output feedback,''
\textit{Automatica}, 41(2), 339--344, 2005.

\bibitem{Michiels04ima}
W. Michiels, S. I. Niculescu, and L. Moreau,
``Using delays and time-varying gains to improve the static
output feedback stabilizability of linear systems: A comparison,''
\textit{IMA Journal of Mathematical Control and Information},
21(4), 393--418, 2004.

\bibitem{boikov05arc}
I. V. Boikov, ``The Brockett stabilization problem,''
\textit{Automation and Remote Control}, 66(5), 746--751, 2005.

\bibitem{Insperger08ifac}
T. Insperger and G. St{\'e}p{\'a}n,
``Brockett problem for systems with feedback delay,''
\textit{IFAC Proceedings Volumes}, 41(2), 11491--11496, 2008.

\bibitem{Insperger06tcst}
T. Insperger, ``Act-and-wait concept for continuous-time
control systems with feedback delay,''
\textit{IEEE Transactions on Control Systems Technology},
14(5), 974--977, 2006.

\bibitem{Gomez19auto}
M. A. Gomez, W. Michiels, and S. Mondi\'e,
``Design of delay-based output-feedback controllers optimizing
a quadratic cost function via the delay Lyapunov matrix,''
\textit{Automatica}, 107, 146--153, 2019.

\bibitem{nicule04tac}
S. I. Niculescu and W. Michiels,
``Stabilizing a chain of integrators using multiple delays,''
\textit{IEEE Transactions on Automatic Control},
49(5), 802--807, 2004.

\bibitem{kharitonov05tac}
V. L. Kharitonov, S. I. Niculescu, J. Moreno, and W. Michiels,
``Static output feedback stabilization: Necessary conditions for multiple delay controllers,''
\textit{IEEE Transactions on Automatic Control},
50(1), 82--86, 2005.

\bibitem{michiels2007}
W. Michiels and S. I. Niculescu,
\textit{Stability and Stabilization of Time-Delay Systems: An Eigenvalue-Based Approach}.
Society for Industrial and Applied Mathematics, Philadelphia, 2007.

\bibitem{karafyllis08scl}
I. Karafyllis,
``Robust global stabilization by means of discrete-delay output feedback,''
\textit{Systems \& Control Letters}, 57(12), 987--995, 2008.

\bibitem{emiliabook}
E. Fridman,
\textit{Introduction to Time-Delay Systems: Analysis and Control}.
Birkh\"{a}user, Basel, 2014.

\bibitem{fridman25delaycontrol}
E. Fridman and A. Selivanov,
``Using delay for control,''
\textit{Annual Review of Control, Robotics, and Autonomous Systems},
8(1), 101--126, 2025.


\bibitem{emilia16auto}
E. Fridman and L. Shaikhet,
``Delay-induced stability of vector second-order systems
via simple Lyapunov functionals,''
\textit{Automatica}, 74, 288--296, 2016.

\bibitem{emilia17auto}
E. Fridman and L. Shaikhet,
``Stabilization by using artificial delays: An LMI approach,''
\textit{Automatica}, 81, 429--437, 2017.




\bibitem{anton18auto}
A. Selivanov and E. Fridman,
``An improved time-delay implementation of derivative-dependent
feedback,'' \textit{Automatica}, 98, 269--276, 2018.

\bibitem{efimov20automatica}
D. Efimov, E. Fridman, W. Perruquetti, and J.-P. Richard,
``Homogeneity of neutral systems and accelerated stabilization
of a double integrator by measurement of its position,''
\textit{Automatica}, 118, 109023, 2020.

\bibitem{iftime26}
O. V. Iftime and S. I. Niculescu,
``Hyperexponential stabilization of chains of integrators using time-varying delayed output feedback,''
\textit{Systems \& Control Letters}, 216, 106544, 2026.

\bibitem{zhoutac}
B. Zhou, W. Michiels, and J. Chen,
``Fixed-time stabilization of linear delay systems by smooth
periodic delayed feedback,''
\textit{IEEE Transactions on Automatic Control},
67(2), 557--573, 2022.

\bibitem{zhang26tac}
K. K. Zhang, E. Fridman, and P. Wang,
``Fixed-time stabilization of systems with nilpotent matrices
in the presence of large unknown delays,''
\textit{IEEE Transactions on Automatic Control}, provisionally accepted (arXiv preprint arXiv:2607.24595), 2026.

\bibitem{Mazenc26fixed}
F. Mazenc and M. Malisoff,
``Fixed-time stabilization of uncertain continuous-time systems with sampled control,''
\textit{IEEE Transactions on Automatic Control},
71(1), 644--651, 2026.

\bibitem{zhang26auto}
K. K. Zhang and E. Fridman,
``Using artificial delays for prescribed/fixed-time static
output feedback control of linear systems,''
arXiv preprint arXiv:2610.33878, 2026.

\bibitem{Zhou23automatica}
B. Zhou,
``On the relative degree and normal forms of linear systems
by output transformation with applications to tracking,''
\textit{Automatica}, 148, 110800, 2023.

\bibitem{zhou2021tac}
B. Zhou and Y. Shi,
``Prescribed-time stabilization of a class of nonlinear systems by linear time-varying feedback,''
\textit{IEEE Transactions on Automatic Control},
66(12), 6123--6130, 2021.

\bibitem{song2017auto}
Y. Song, Y. Wang, J. Holloway, and M. Krstic,
``Time-varying feedback for regulation of normal-form nonlinear systems in prescribed finite time,''
\textit{Automatica},
83, 243--251, 2017.

\bibitem{ye2023tac}
H. Ye and Y. Song,
``Prescribed-time tracking control of MIMO nonlinear systems with nonvanishing uncertainties,''
\textit{IEEE Transactions on Automatic Control},
68(6), 3664--3671, 2023.

\bibitem{polyakov12}
A. Polyakov,
``Nonlinear feedback design for fixed-time stabilization of linear control systems,''
\textit{IEEE Transactions on Automatic Control},
57(8), 2106--2110, 2012.

\bibitem{davison73tac}
E. J. Davison and S. H. Wang,
``Properties of linear time-invariant multivariable systems
subject to arbitrary output and state feedback,''
\textit{IEEE Transactions on Automatic Control},
18(1), 24--32, 1973.

\bibitem{horn13}
R. A. Horn and C. R. Johnson,
\textit{Matrix Analysis}.
Cambridge University Press, 2012.

\bibitem{karafyllis11tac}
I. Karafyllis and M. Krstic,
``Nonlinear stabilization under sampled and delayed measurements, and with inputs subject to delay and zero-order hold,''
\textit{IEEE Transactions on Automatic Control},
57(5), 1141--1154, 2012.

\bibitem{zhang22tac}
J. Zhang and E. Fridman,
``Constructive robust stabilization by using square wave
dithers: A time-delay approach,''
\textit{IEEE Transactions on Automatic Control},
68(10), 6307--6314, 2023.


\bibitem{qiu93auto}
L. Qiu and E. J. Davison,
``Performance limitations of non-minimum phase systems
in the servomechanism problem,''
\textit{Automatica}, 29(2), 337--349, 1993.

\bibitem{CTMSpole}
D. Tilbury and B. Messner,
``Inverted pendulum: System modeling,''
2011. [Online]. Available:
\url{https://ctms.engin.umich.edu/CTMS/index.php?example=InvertedPendulum&section=SystemModeling}.
Accessed: Aug. 8, 2026.

\end{thebibliography}
\end{document}